\documentclass[12pt,a4paper,fleqn]{article}
\usepackage{tikz}
\usepackage{caption}
\usepackage{lscape}

\usepackage{fullpage}
\usepackage{amsmath}
\usepackage{amsthm}
\usepackage{amssymb}
\usepackage{amscd}
\usepackage{plain}
\usepackage{graphicx}
\usepackage{t1enc,epsfig}
\usepackage[mathscr]{eucal}
\usepackage{epstopdf}
\usepackage[german,english]{babel}
\usepackage{stmaryrd}

\usepackage{xcolor}

\usepackage{wrapfig}
\usepackage{float}

\usepackage{ulem}
\usepackage{cancel}

\newtheorem{theorem}{Theorem}[section]

\newtheorem{corollary}[theorem]{Corollary}
\newtheorem{definition}[theorem]{Definition}
\newtheorem{example}[theorem]{Example}
\newtheorem{lemma}[theorem]{Lemma}
\newtheorem{proposition}[theorem]{Proposition}
\newtheorem{remark}[theorem]{Remark}
\newtheorem{conjecture}[theorem]{Conjecture}
\numberwithin{equation}{section}

\begin{document}

\newcommand {\beq}{\begin{equation}}
\newcommand {\eeq}{\end{equation}}
\newcommand{\bthm}{\begin{theorem}}
\newcommand{\ethm}{\end{theorem}}
\newcommand{\bd}{\begin{definition}}
\newcommand{\ed}{\end{definition}}
\newcommand{\bs}{\begin{proposition}}
\newcommand{\es}{\end{proposition}}
\newcommand{\bp}{\begin{proof}}
\newcommand{\ep}{\end{proof}}
\newcommand{\br}{\begin{remark}}
\newcommand{\er}{\end{remark}}
\newcommand{\bex}{\begin{example}}
\newcommand{\eex}{\end{example}}
\newcommand{\bc}{\begin{corollary}}
\newcommand{\ec}{\end{corollary}}
\newcommand{\bl}{\begin{lemma}}
\newcommand{\el}{\end{lemma}}
\newcommand{\bfj}{\begin{conjecture}}
\newcommand{\ej}{\end{conjecture}}

\newcommand{\red}[1]{{\color{red}#1}}

\def\cl{\centerline}
\def\diy{\displaystyle}
\def\cre{\color{red}} \def\cbl{\color{blue}} \def\cte{\color{teal}}

\cl{\large{\bf A large deviation principle for birth-death processes}}\vskip .3cm
\cl{\large{\bf with a linear rate of downward jumps}}\vskip .5cm

\cl{\bf{N. Vvedenskaya$^{1)}$, A. Logachov$^{2,3)}$, Y. Suhov$^{4,5)}$,
 A.Yambartsev$^{2)}$}}\vskip .5cm

\footnote{
	
	\noindent$^{1)}$ {\footnotesize Institure for Information Transmission Problems, RAS, Moscow, RF}
	
	\noindent$^{2)}$ {\footnotesize IME Universidade de Sao Paulo, Brazil}
	
	\noindent$^{3)}$ {\footnotesize Sobolev Institute of Mathematics, RAS Siberian Branch, Novosibirsk, RF}
	
	\noindent$^{4)}$ {\footnotesize Department of Mathematics, Penn State University, University Park,
		State College, PA, USA}
	
	\noindent$^{5)}$ {\footnotesize DPMMS, University of Cambridge,  and St John's College, Cambridge, UK}
}


{\large\bf Abstract.} Birth-death processes form a natural class where ideas and results
on large deviations can be tested. In this paper, we derive a large deviation principle
under an assumption that the rate of a jump down (death) is growing asymptotically
linearly with the population size, while the rate of a jump up (birth) is growing sub-linearly.
 We establish a large deviation principle under  various forms of scaling of the underlying
process and the corresponding normalization of the logarithm of the large deviation probabilities.
The results show interesting features of dependence of
the large deviation functional upon the parameters of the process and the forms of scaling
and normalization.

\section{Introduction and definitions}\label{Sect1} \def\lam{\lambda} \def\veps{\varepsilon}

A birth-death process is a continuous-time Markov process with states $x\in\mathbb{Z}^+
:=\{0,1,2,\ldots\}$ (representing the population size) and with transitions occurring between
neighboring states. The class of birth-death processes exhibits a remarkable balance between
simplicity, allowing for analytical solutions, and complexity, showcasing a diverse range of
interesting phenomena. Its versatility is
accentuated by the possibility of exploring various jump rates, drawing attention from multiple
research areas. Furthermore, birth-death processes find applications across diverse fields, such
as information theory (involving encoding and storage of information, see \cite{SS1}), population
biology, genetics, ecology, (reviewed in    \cite{BioRev1}, \cite{BioRev2}), chemistry (modeling
growth and extinction in systems with multiple components,  see \cite{SS2},  \cite{MSSZ}),
economics (modeling competitive production and pricing, \cite{MPY}, \cite{VSB}) and queueing
system theory (explored, for example,  in \cite{Rob}).

In particular, birth-death processes are instrumental in exploring various aspects
of the large deviation theory, which is the focus of this paper. Apart from that, we mention
connections between birth-death processes and orthogonal polynomials, as detailed in \cite{Val1}
and \cite{Doorn}.

In this paper we work under the assumption that the rate $\lam (x)$ of jump $x\to x+1$ and
the rate $\mu (x)$ of jump $x\to x-1$ obey the condition \eqref{17.12.3}: $\mu (x)$ grows
with $x$ asymptotically linearly, while $\lam (x)$ grows asymptotically sub-linearly.
This assumption ensures positive-recurrence of the process (cf. \cite{Karl5}). Such processes
find an application in population dynamics \cite{Kar}, \cite{Kend1}; they are also relevant in
models of a market interaction between ask-bid sides of a limit order book \cite{MPY},
energy-efficient schemes for cloud resources
\cite{Shi} and scenarios with an increasing number of available servers in stations \cite{Rob}.

\vskip .2cm







Let us provide formal definitions. We consider a continuous-time Markov process $\xi (t)$, $t\geq 0$,
on the state space $\mathbb{Z}^+$, starting at point $0$. The process dynamics is
as follows. There are given two functions, $\lam :\,\mathbb{Z}^+\to (0,\infty )$ giving the rate of upward jumps, and
$\mu\,:\,\mathbb{Z}^+\to [0,\infty )$ giving the rate of downward jumps, with $\mu (0)=0$ and $\mu (x)>0$ for $x\geq 1$.
We set $\eta =\lam +\mu$ for the combined jump rate. Given that $\xi (t)=x$ for some $t\geq 0$ and $x\in\mathbb{Z}^+$,
the value of the process remains unchanged for an exponentially distributed random time $\tau_x$ of rate $\eta (x)$.
At time $t+\tau_x$ the process jumps to either $x+1$ or $x-1$ with the following probabilities
\beq \label{17.12.1}
{\bf P}(\xi (t+\tau_x)=x+1)=\dfrac{\lambda(x)}{\eta(x)},\ \ \ {\bf P}(
 \xi(t+\tau_x)=x-1)=\dfrac{\mu(x)}{\eta(x)}.\eeq
For the case where $x=0$, the only feasible transition is to site $1$.
The key assumption is that there exist constants $P,Q>0$ and $l\in[0,1)$ such that
\beq \label{17.12.3}
\lim_{x\to\infty}\dfrac{\lambda(x)}{x^l}=P, \ \ \ \lim_{x\to\infty}\dfrac{\mu(x)}{x}=Q.\eeq

We focus on the large deviation principle (LDP) for the family of processes
\beq\label{17.12.4}
\xi_T(t)=\dfrac{\xi(tT)}{\varphi(T)}, \ 0\leq t\leq 1,\eeq
for a subexponential (\ref{17.12.5}), exponential (\ref{11.02.7}) or superexponential (\ref{17.12.6})
growth of value $\varphi(T)$. Here $T>0$ is a time-scaling parameter, and
$\varphi :(0,\infty )\to (0,\infty)$ is a Lebesgue-measurable function
referred to as a scaling function. We assume that
\beq\label{17.19.19}
\lim\limits_{T\rightarrow\infty}\varphi(T)=\infty .\eeq

The space where we will establish the large deviation principle is $\mathbb{L}=\mathbb{L}_1[0,1]$,
with the standard metric
$$\rho(f,g)=\int\limits_0^1|f(t)-g(t)|dt, \ f,g\in\mathbb{L}.$$
Let $\mathfrak{B}=\mathfrak{B}_{(\mathbb{L},\rho)}$ denote the Borel $\sigma$\,-\,algebra
in $(\mathbb{L},\rho)$; for a set $\mathbb{B}\in\mathfrak{B}$, ${\rm{cl}}(\mathbb{B})$
and ${\rm{int}}(\mathbb{B})$ stand for the closure and the interior of $\mathbb{B}$, respectively.

Recall the notions and definitions we need (see for more details \cite{DZ}--\cite{Feng}, \cite{Stro},
\cite{Var1}, \cite{Var2}).
In Definitions \ref{d1.1} and \ref{d1.2} below we attempt to cover a variety of situations occurring
in the context of the current paper. In these definitions we use a Lebesgue-measurable function $\psi :(0,\infty )\to (0,\infty )$ satisfying $$\lim\limits_{T\to\infty}\psi (T)=\infty,$$ and a $\mathfrak{B}$-measurable
functional $I \,:\, \mathbb{G}\rightarrow [0,\infty ]$ where $\mathbb{G}\subseteq\mathbb{L}$
and $\mathbb{G}\in\mathfrak{B}$.  Given $\mathbb{A}\subseteq\mathbb{G}$
such that $\mathbb{A}\in
\mathfrak{B}$, we set $$I(\mathbb{A}) = \inf\limits_{y \in\mathbb{A}} I(y),$$ with $I(\varnothing)=\infty$.
Furthermore, $\psi$ is referred to as a normalizing function and $I$ as a large deviation (LD)  rate functional.
\vskip .5cm

\begin{definition} \label{d1.1}
Let $\mathbb{G}\subseteq\mathbb{L}$ and $\mathbb{G}\in\mathfrak{B}$. Let a family of random
processes $\xi_T(\,\cdot\,)$, $T>0$, be
defined as in \eqref{17.12.4} for some scaling function $\varphi$.
We say that this family satisfies an $(\mathbb{G},\mathbb{L},\rho)$-local large deviation principle
$\big((\mathbb{G},\mathbb{L},\rho)$-{\rm{LLDP}}$\big)$ with an {\rm{LD}} functional
$I\,:\mathbb{G}\to [0,\infty ]$  and the normalizing function $\psi$ if \ for all \ $f \in\mathbb{G}$
$$\begin{array}{l}
\lim\limits_{\veps\rightarrow 0}\limsup\limits_{T\rightarrow \infty}\dfrac{1}{\psi(T)}
\ln\mathbf{P}(\xi_T(\,\cdot\,)\in\mathbb{U}_\veps(f))\\
\qquad =\lim\limits_{\veps\rightarrow 0}\liminf\limits_{T\rightarrow \infty}\dfrac{1}{\psi(T)}
\ln\mathbf{P}(\xi_T(\,\cdot\,)\in\mathbb{U}_\veps(f))=-I(f),\end{array}$$
where
$$\mathbb{U}_\veps(f)=\{g\in \mathbb{L}: \ \rho(f,g)<\veps\}.$$
\end{definition} \vskip .5cm

\begin{definition} \label{d1.2} Let $\xi_T(\,\cdot\,)$, $T>0$, be family of random processes  defined
as in \eqref{17.12.4} for some scaling function $\varphi$. We say that this family satisfies an
$(\mathbb{L},\rho )$-{\rm{LDP}} with a normalizing function
$\psi$ and an {\rm{LD}} functional $I\,:\;\mathbb{L}\to (0,\infty ]$ if, whenever set $\mathbb{B}
\subseteq\mathbb{L}$ and $\mathbb{B}\in \mathfrak{B}$, we have that
$$\begin{array}{c} \limsup\limits_{T \rightarrow \infty} \dfrac{1}{\psi(T)} \ln \mathbf{P}(\, \xi_T(\,\cdot\,)
\in\mathbb{B} \,) \leq - I({\rm{cl}}(\mathbb{B})),\\
\liminf\limits_{T \rightarrow \infty} \dfrac{1}{\psi(T)} \ln \mathbf{P}(\, \xi_T(\,\cdot\,)
\in\mathbb{B}\,)\geq -I({\rm{int}}(\mathbb{B})).\end{array}$$
\end{definition} \vskip .5cm

\begin{definition} \label{d1.3} Let $\xi_T(\,\cdot\,)$, $T>0$, be family of random processes  defined
as in \eqref{17.12.4} for some scaling function $\varphi$. We say that this family is exponentially
tight {\rm (ET)} on $(\mathbb{L},\rho )$ with a normalizing function $\psi$  if for any $C>0$ there
exists a compact set  $\mathbb{K}_C\subseteq\mathbb{L}$ such that
$$\limsup\limits_{T \rightarrow \infty} \dfrac{1}{\psi(T)} \ln \mathbf{P}(\, \xi_T(\,\cdot\,)
\not\in\mathbb{K}_C \,) \leq - C.$$
\end{definition} \vskip .5cm

If a family $\xi_T(\,\cdot\,)$, $T>0$, and a functional $I$ satisfy Definitions \ref{d1.2} and \ref{d1.3}
(in particular, family $\xi_T(\,\cdot\,)$, $T>0$, is ET) then for all $c\geq 0$  the set
$\{ f \in \mathbb{L}\,:\, I(f) \leq c \}$ is a compact in $(\mathbb{L},\rho)$. In this case,
one says that $I$ is a ``good rate functional''  (cf. \cite[section 1.2]{DZ}, \cite[section 2.2]{Feng}).
In this paper the ET property is established in Lemma \ref{l4.5}. It is known (see, for example,
\cite{Pukh}) that if the trajectories of random processes $\xi_T(\,\cdot\,)$ belong to a Polish space
 then the ET property is a necessary condition for the goodness of functional $I$.
Note that this holds true in our setting.

An earlier  work \cite{VLSY}
established an LLDP  for a family of processes \eqref{17.12.4} with the scaling function $\varphi(T)=T$, while
paper \cite{VLSY1} did it for the case of subexponential asymptotics of $\varphi(T)$, when
\beq\label{17.12.5}
\lim\limits_{T\rightarrow\infty}\dfrac{\ln\varphi(T)}{T}=0.
\eeq
In this latter case, the family \eqref{17.12.4} is not ET (we discuss this in Section~\ref{Discussion}).
Consequently, the LDP is not available in the whole of $(\mathbb{L},\rho)$. 

In the present paper we consider two complementary conditions:

1) there exists a constant $k\in (0,\infty )$ such that
\beq\label{11.02.7}
\lim\limits_{T\rightarrow\infty}\dfrac{\ln\varphi(T)}{T}=k,
\eeq

and 2) the limit
\beq\label{17.12.6}
\lim\limits_{T\rightarrow\infty}\dfrac{\ln\varphi(T)}{T}=\infty.
\eeq
The form of the LD functional depends on which condition is assumed,
\eqref{11.02.7} or \eqref{17.12.6}, cf. Section 2, Theorems  \ref{th2.1} and \ref{th2.2}. An emerging question
is: Why do scalings \eqref{11.02.7} or \eqref{17.12.6} lead to the large deviation principle, while scaling
\eqref{17.12.5} does not? We will explain this in Section~\ref{Discussion}.
\vskip .2cm

Let us discuss what is currently known outside condition \eqref{17.12.3}; cf. \cite {VLSY1}. Suppose that
$\lam (x)\sim Px^l$ and $\mu (x)\sim Qx^m$ where $0\leq l<m$. If $m\in (0,1)$ then three cases emerge,
depending on a condition upon scaling function $\varphi$, and the form of the rate functional is different in
each of these cases.
If we assume that $m>1$ then only an LLDP will take place, so the three cases will be
reduced to one.  Also, from  \cite{VLSY1} it follows that an LLDP holds true
when rates $\lam (x)$ and $\mu (x)$ are regularly varying functions.
Separately, notice the case where $\lam (x)=P$, $\mu (x)=Q$ where $P$ and $Q$ are positive
constants. Here, process $\xi (t)$ is compound Poisson, for which the LD asymptotics are well-known
\cite{BorMog2}, \cite{Lynch}, \cite{Mog2}.

\vskip .2cm

This paper contains four sections. In Section 2 we state our main result, Theorem  \ref{th2.2},
and a triple of auxiliary assertions (Lemmas 2.\ref{l2.3} -- \ref{l2.5}).
Section 3 is dedicated to the derivation of Theorem \ref{th2.2} from Lemmas 2.\ref{l2.3} -- \ref{l2.5}
and the proof of these lemmas. Section~\ref{Discussion}
contains a discussion of obtained results. Finally, in Section 5 we prove some additional technical assertions
(Lemmas \ref{l4.2} -- \ref{l4.5}) used in the proof or an interpretation of the obtained results.

\vskip .2cm

{\it A commemorative note.} It is with great sadness and sorrow that the rest of the authors report of the loss of our remarkable collaborator and friend Nikita Vvedenskaya (1930-2022). Until her last days she actively worked on this project, and her contribution was essential and irreplaceable. We will miss her dearly.

\section{Notation, the main result}

We denote by $\mathbb{V}=\mathbb{V}[0,1]$  the set of non-negative measurable  functions
$f:\,[0,1]\mapsto [0,\infty)$ of a finite variation.
Given $f\in\mathbb{V}$, let $\mathrm{Var}\,f$ be the total variation of $f$.

Next, $\mathbb{C}=\mathbb{C}[0,1]$ is the space of continuous functions on $[0,1]$. From now on we
let $\mathbb{G}$ be the set of functions $f\in\mathbb{C}$ such that $f(0)=0$ and $f(t)> 0$  for $t>0$.

The following result follows from \cite{VLSY1}. \vskip .5cm

\begin{theorem} \label{th2.1}
Assume conditions \eqref{17.12.3} and \eqref{17.12.5}. Then the  family $\xi_T(\,\cdot\, )$, $T>0$,
defined as in \eqref{17.12.4} satisfies an $(\mathbb{G},\mathbb{L},\rho)$-{\rm{LLDP}}
with the normalizing
function $\psi(T)=T\varphi(T)$ and the {\rm{LD}} functional
$$I(f)=Q\int_0^1f(t)dt.$$ \end{theorem} \vskip .5cm

Given $f\in\mathbb{V}$, we use the decomposition into monotone increasing and decreasing
components
\beq\label{fpm} f(t)=f^+(t)-f^-(t), \  f^+(0)=f(0), \  f^-(0)=0, \ \hbox{with} \
\mathrm{Var}\,f=\mathrm{Var}\,f^++\mathrm{Var}\,f^-.\eeq
Such a decomposition is unique (cf. \cite[Ch 1, \S 4]{Rise}).

Denote by $\mathbb{D}=\mathbb{D}[0,1]$ the space of c\`adl\`ag functions on $[0,1)$ with left limits at $t=1$.
Observe that for every $f\in\mathbb{V}$ there exists a function $f_\mathbb{D}\in\mathbb{D}$
such that $\rho(f,f_\mathbb{D})=0$.

The main result of this paper is \vskip .5cm

\begin{theorem} \label{th2.2}
Assume condition \eqref{17.12.3}.

{\rm{1)}} Under condition \eqref{11.02.7} the family  $\xi_T(\,\cdot\,)$, $T>0$,
defined as in \eqref{17.12.4} satisfies an
$(\mathbb{L},\rho)$-{\rm{LDP}} with the normalizing function $\psi(T)=\varphi(T)\ln\varphi(T)$
and the good {\rm{LD}} functional $I\,:\,\mathbb{L}\to [0,\infty ]$ where
$$I(f)=\left\{ \begin{aligned}
& \frac{Q}{k}\diy\int_0^1f(s)ds+(1-l)f^+_\mathbb{D}(1), \ f\in \mathbb{V},\\
& \infty,   \ f\notin \mathbb{V}.\end{aligned} \right.$$

{\rm{2)}} Under condition \eqref{17.12.6} the family  $\xi_T(\,\cdot\,)$, $T>0$,
defined as in \eqref{17.12.4} satisfies an $(\mathbb{L},\rho)$-{\rm{LDP}} with
the normalizing function $\psi(T)=\varphi(T)\ln\varphi(T)$  and the
good {\rm{LD}} functional $I\,:\,\mathbb{L}\to [0,\infty ]$ where
$$I(f)=\left\{ \begin{aligned}
& (1-l)f^+_\mathbb{D}(1), \ f\in \mathbb{V},\\
&\infty,   \ f\notin \mathbb{V}.\end{aligned} \right.$$
\end{theorem} \vskip .5cm

Before we pass to the proof, let us make some comments.
Note that the LDP in the space of the right-continuous functions with the Skorokhod metric is not obtained
since the set of functions with the total variation bounded by a constant is non-compact in this space. On the other hand, it seems that the results of this paper will hold for the space of functions without second-kind discontinuities equipped with the Borovkov metric (cf. \cite{BorMog3}--\cite{BorMog5}).
It is also worth mentioning that, in contrast with the classical results, in our case the LD functional $I(f)$
does not contain the integral of the convex function of the derivative of the absolutely
continuous component of the function $f$.
\vskip .2cm

The proof  of Theorem \ref{th2.2} uses auxiliary assertions; see Lemmas \ref{l2.3} -- \ref{l2.5} below.
Let us introduce some additional notions.
Given $T>0$, denote by $\mathbb{X}_T$ the set of right-continuous functions $u\,:\,[0,T]\to\mathbb{Z}^+$
with $u(0)=0$, having a finite number of jumps $n (u)$, where every jump has size $\pm 1$. This
gives the set of trajectories for the birth-death process $\xi (t)$, $t\in [0,T]$. We speak below of measures
on $(\mathbb{X}_T,\mathfrak{X}_T)$, where $\mathfrak{X}_T$ is a standard Borel
$\sigma$\,-\,algebra in $\mathbb{X}_T$.

Next, consider a continuous-time Markov process $\zeta(t)$, $t\in[0,T]$, on the state-space
$\mathbb{Z}$, with the full jump rate $1$, jump size $\pm 1$, and probabilities of jumps $1/2$.
There is a positive probability that this process lives in $\mathbb{X}_T$. In Lemma  2.\ref{l2.3}
and later we refer to the two processes as $\xi$ and $\zeta$.

\begin{lemma}\label{l2.3} \rm{(cf. \cite{MPY}, \cite{VLSY})} The distribution of the random
process $\xi$ on $\mathbb{X}_T$ is absolutely continuous with respect
to that of a process $\zeta$. The corresponding Radon--Nikodym density
$\mathbf{p}=\mathbf{p}_T$ on $\mathbb{X}_T$ has the form:
\beq\label{17.12.7}
\mathbf{p}(u)=\left\{ \begin{array}{ll}2^{n(u)}\Bigl(\prod\limits_{i=1}^{n(u)}
e^{-(\eta(u(t_{i-1}))-1)\tau_{i}}\nu(u(t_{i-1}),u(t_i))\Bigr)&\;\\
\qquad \times e^{-(\eta(u(t_{n(u)})-1))(T-t_{n(u)})},&\text{if } \;n(u)\geq 1,\\
e^{-(\eta(0)-1)T},&\text{if } \;  n(u)= 0,\end{array} \right.\eeq
where $\eta (x) = \lam (x) + \mu (x)$, $x\in\mathbb{Z}^+$; cf. \eqref{17.12.1}. Here we suppose
that function $u\in\mathbb{X}_T$  has jumps at time-points  \ $0<t_1<...<t_{n(u)}<T$ and set
$\tau_i=t_i-t_{i-1}$, with $t_0=0$. Further, the value $\nu (u(t_{i-1}),u (t_i))$ is given by
$$
\nu (u(t_{i-1}),u (t_i))= \left\{ \begin{aligned}\lambda(u (t_{i-1})), &
\mbox{ if }\; u (t_i)-u (t_{i-1})=1;\\
\mu(u (t_{i-1})), & \mbox{ if }\;  u (t_i)-u (t_{i-1})= -1. \end{aligned} \right.$$
\end{lemma} \vskip .5cm

Let $N_T(\zeta)$ be the number of jumps in process $\zeta(t)$ on the interval $[0,T]$.
The claim of  Lemma~2.\ref{l2.3} is equivalent to the fact that for any measurable set
$\mathbb{H}\subseteq\mathbb{X}_T$
\beq\label{17.12.8}
\mathbf{P}(\xi\in\mathbb{H})=e^T \mathbf{E}\bigl[ e^{-A_T(\zeta)}e^{B_T(\zeta)+N_T(\zeta)\ln2}
{\mathbf 1}(\zeta\in\mathbb{H})\bigr].\eeq
Here
\beq\label{17.12.9}\begin{aligned} A_T(\zeta):=& \diy\int_0^T \eta(\zeta(t))dt\\
 = & \begin{cases}\sum\limits_{i=1}^{N_T(\zeta)}\eta(\zeta(t_{i-1}))\tau_{i}+
\eta(\zeta(t_{N_T(\zeta)}))(T-t_{N_T(\zeta)}), &
\mbox{if } N_T(\zeta)\geq1,\\
\eta(0)T, & \mbox{if }\; N_T(\zeta)=0,\end{cases}\end{aligned}\eeq
and
\beq\label{17.12.10}
B_T(\zeta):= \begin{cases}\sum\limits_{i=1}^{N_T(\zeta)}\ln(\nu(\zeta(t_{i-1}),\zeta(t_i))), &
\mbox{if }\; N_T(\zeta)\geq1;\\
0, & \mbox{if }\;  N_T(\zeta)=0. \end{cases}\eeq
Symbols $\mathbf{1}(\;\cdot\;)$ and $\mathbf{1}[\;\cdot\;]$ stand for the indicators of events in
$\sigma$\,-\,algebra $\mathfrak{B}$.

Representation \eqref{17.12.8} is used in the analysis of the value $\ln \mathbf{P}(\xi_{T}(\cdot)\in
\mathbb{U}_\veps(f))$.  We set
$$\zeta_T(t):=\dfrac{\zeta(tT)}{\varphi(T)}, \ \ t\in[0,1].$$
In what follows, we write $\xi_T, \zeta_T$ instead of $\xi_T(\cdot), \zeta_T(\cdot)$  and $A_T$, $B_T$,  $N_T$
instead of $A_T(\zeta)$, $B_T(\zeta)$, $N_T(\zeta)$.

The proof of Theorem \ref{th2.2} is based on the analysis  of $\mathbf{p}_T$. This is a common
method in the LD theory, particularly, in a specification of an LD functional. Namely, we analyze the
Radon-Nikodym density  $\mathbf{p}_T$ on the event $\big\{\zeta_T\in\mathbb{U}_\veps(f)\big\}$
and prove an $(\mathbb{V},\mathbb{L},\rho )$-LLDP by using the independence of increments in process $\zeta_T$,
together with the Stirling formula and properties of the functional space $(\mathbb{L},\rho )$. 
Cf. Lemmas \ref{l2.4} and \ref{l2.5} below and their proof in Section 5. Next, we prove that the family $\xi_T$
is ET (cf. Lemma \ref{l4.5} in Section 5).  Then, by using a standard implication LLDP plus ET $\Rightarrow$ LDP
(cf. \cite[Lemma 4.1.23]{DZ}, \cite{Pukh}), we obtain an $(\mathbb{L},\rho)$-LDP for processes $\xi_T$.


\begin{lemma}\label{l2.4} Assume condition \eqref{17.12.3} and one of conditions \eqref{11.02.7}
or \eqref{17.12.6}. Then for all $f\in \mathbb{V}$ with $\rho(f,0)> 0$,
$$\lim\limits_{\veps\rightarrow 0} \limsup\limits_{T\rightarrow\infty}
\dfrac{\ln\mathbf{E}\bigl[e^{B_T+N_T\ln2}\mathbf{1}(\zeta_T\in
\mathbb{U}_\veps(f))\bigr]}{\varphi(T)\ln\varphi (T)}\leq (1-l)f^+_\mathbb{D}(1).$$ \end{lemma}\vskip .5cm

\begin{lemma} \label{l2.5} Assume condition \eqref{17.12.3} and one of conditions \eqref{11.02.7}
or \eqref{17.12.6}. Then for all $f\in \mathbb{V}$ with $\rho(f,0)> 0$,
$$\lim\limits_{\veps\rightarrow 0} \liminf\limits_{T\rightarrow\infty}
\dfrac{\ln\mathbf{E}\bigl[e^{B_T+N_T\ln2}{\mathbf 1}(\zeta_T\in\mathbb{U}_\veps(f))\bigr] }{\varphi(T)
\ln\varphi(T)}\geq (1-l)f^+_\mathbb{D}(1).$$\end{lemma}

\section{Proof of Theorem \ref{th2.2} and Lemmas \ref{l2.4},  \ref{l2.5}}

\smallskip

\quad\; Symbol $\Box$ marks the end of a proof.
\vskip .5cm

{\sc Proof of Theorem \ref{th2.2}.} First, consider the case where $\rho(f,0)=0$. Obviously,
\beq \label{12.10.2}
\lim\limits_{\veps\rightarrow 0} \limsup\limits_{T\rightarrow\infty}
\dfrac{1}{\varphi(T)\ln\varphi(T)}\ln\mathbf{P}\bigl(\xi_T\in  \mathbb{U}_\veps(f)\bigr)\leq 0 = I(f). \eeq
It is easy to see that
\beq\label{12.10.1}\begin{aligned}
\lim\limits_{\veps\rightarrow 0} \liminf\limits_{T\rightarrow\infty}
& \dfrac{1}{\varphi(T)\ln\varphi(T)}\ln\mathbf{P}\bigl(\xi_T\in  \mathbb{U}_\veps(f)\bigr) \\
& \geq \lim\limits_{\veps\rightarrow 0} \liminf\limits_{T\rightarrow\infty}
\dfrac{1}{\varphi(T)\ln\varphi(T)}\ln\mathbf{P}\Bigl(\sup\limits_{t\in[0,1]}\xi_T(t)=0\Bigr)\\
& =\lim\limits_{\veps\rightarrow 0} \liminf\limits_{T\rightarrow\infty}\dfrac{1}{\varphi(T)
\ln\varphi(T)}\ln e^{-\lambda(0)T}=0.\end{aligned}
\eeq
Now, suppose that $\rho(f,0)>0$. We start with evaluating $A_T$. As follows from \eqref{17.12.9},
\beq\label{11.02.51}
A_T:=\int_0^T \eta(\zeta(t))dt=T\int_0^1 \eta(\varphi(T)\zeta_T(s))ds. \eeq
Condition \eqref{17.12.3} implies that, for any given $\varepsilon, \gamma\in(0,1)$, for $T$ large enough
\beq\label{17.12.16}
Q(1-\gamma)\leq\dfrac{\eta\bigl( \varphi(T)(\zeta_T(s)\vee\veps) \bigr)}
{\varphi(T)(\zeta_T(s)\vee\veps)} \leq Q(1+\gamma).\eeq
Here and below, $a\vee b=\max(a,b)$. Furthermore, the values $\varepsilon$ and $\gamma$ will tend to zero.

Let us upper-bound the integral in \eqref{11.02.51}. Suppose $\zeta _T\in  \mathbb{U}_\veps(f)$.
Then the right
bound in \eqref{17.12.16} implies that for any $\veps,\gamma\in(0,1)$, if $T$ is large enough,
we have the inequalities
\beq\label{11.02.8}\begin{aligned}\diy\int_0^1 \eta(\varphi(T)\zeta_T(s))ds
 & \leq\int_0^1 \eta(\varphi(T)(\zeta_T(s)\vee\veps))ds\\
& \leq\varphi(T)(1+\gamma)Q\diy\int_0^1(\zeta_T(s)\vee\veps)ds\\
& \leq\varphi(T)(1+\gamma)Q\diy\int_0^1(|\zeta_T(s)-f(s)|+f(s)+\veps)ds\\
& \leq\varphi(T)(1+\gamma)Q\diy\int_0^1f(s)ds +2\varphi(T)(1+\gamma)Q\veps.\end{aligned}\eeq

Next, consider a lower bound for the integral in \eqref{11.02.51}. Due to an asymptotic character of
condition (\ref{17.12.3}), we need some caution when dealing with the regions where the scaled process
approaches level zero. Set:
\beq\label{02.08.51}
H:=\{t\in[0,1]:f(t)>0\}, \ \ \ H_\veps:=\{t\in[0,1]:f(t)\geq\veps+\sqrt{\veps}\},
\eeq
\beq\label{02.08.52}
G_\veps:=\{t\in[0,1]:\zeta_T(t)<\veps,f(t)\geq\veps+\sqrt{\veps}\}.
\eeq
If $\zeta_T\in \mathbb{U}_\veps(f)$, the left-hand bound in \eqref{17.12.16} implies that, once more, for any given small $\veps$ and $\gamma$ within the interval $(0,1)$, and with a sufficiently large value of $T$, we have the following:
\beq \label{11.02.9}\begin{array}{l}
\diy\int_0^1 \eta(\varphi(T)\zeta_T(s))ds\geq \varphi(T)(1-\gamma)Q
\int_{H_\veps\setminus G_\veps}\zeta_T(s)ds\\
\quad\geq\varphi(T)(1-\gamma)Q\diy\int_{H_\veps\setminus G_\veps}f(s)ds-
\varphi(T)(1-\gamma)Q\diy\int_{H_\veps\setminus G_\veps}|\zeta_T(s)-f(s)|ds\\
\quad\geq \varphi(T)(1-\gamma)Q\diy\int_{H_\veps\setminus G_\veps}f(s)ds-
\varphi(T)(1-\gamma)Q\veps .\end{array}\eeq

If $\zeta_T\in  \mathbb{U}_\veps(f)$, the Lebesgue measure of the set $G_\varepsilon$
defined by (\ref{02.08.52}) has the following upper bound. Since $f(s)-\zeta_T(s) \ge \sqrt{\varepsilon}$
for all $s\in G_\varepsilon$ we have
\beq\label{11.02.10}
L(G_\veps)=\int_{G_\veps}ds\leq \int_0^1\frac{|\zeta_T(s)-f(s)|}{\sqrt{\varepsilon}}ds =\frac{\rho(\zeta_T,f)}{\sqrt{\varepsilon}} \leq\sqrt{\veps}.\eeq

By virtue of \eqref{17.12.8}, \eqref{11.02.51} and \eqref{11.02.8}, \eqref{11.02.9}, we obtain
that for $T$ large enough
\beq\label{17.12.18}\begin{aligned}
&\exp\bigg\{T-T\varphi(T)(1-\gamma)Q\int_{H_\veps\setminus G_\veps}f(s)ds+
T\varphi(T)(1-\gamma)Q\veps\bigg\} \\
&\qquad\times
\mathbf{E}\bigl[e^{B_T+N_T\ln2}\mathbf{1}(\zeta_T\in\mathbb{U}_\veps(f))\bigr]
\phantom{\int_0^1}\\
&\quad \geq \mathbf{P}\bigl(\xi_T(\cdot)\in\mathbb{U}_\veps(f)\bigr)\phantom{\int_0^1}\\
&\quad \geq
\exp\bigg\{T-T\varphi(T)(1+\gamma)Q\int_0^1f(s)ds -2T\varphi(T)(1+\gamma)Q\veps\bigg\} \\
&\qquad \times
\mathbf{E}\bigl[e^{B_T+N_T\ln2}\mathbf{1}(\zeta_T\in  \mathbb{U}_\veps(f))\bigr].
\phantom{\int_0^1}\end{aligned}\eeq
The bounds \eqref{17.12.18} conclude an initial part of the proof of Theorem \ref{th2.2}. Subsequent parts
establish assertions 1) and 2) based on \eqref{17.12.18}, while assuming conditions \eqref{11.02.7} and
\eqref{17.12.6}, respectively.
\vskip .2cm

First, assume condition \eqref{11.02.7}. According to the upper bound in \eqref{17.12.18}, for any
$\veps,\gamma\in(0,1)$
\beq\label{11.02.11}
\begin{aligned}
\limsup\limits_{T\rightarrow\infty}
\dfrac{1}{\varphi(T)\ln\varphi(T)} & \ln\mathbf{P}\bigl(\xi_T(\cdot)\in  \mathbb{U}_\veps(f)\bigr) \\
\leq -\dfrac{Q(1-\gamma)}{k} & \diy\int_{H_\veps\setminus G_\veps}f(s)ds
+\dfrac{Q(1-\gamma)}{k}\veps \\
+\limsup\limits_{T\rightarrow\infty} &
\dfrac{1}{\varphi(T)\ln\varphi(T)}\ln\mathbf{E}\bigl[e^{B_T+N_T\ln2}
\mathbf{1}(\zeta_T\in  \mathbb{U}_\veps(f))\bigr].
\end{aligned}\eeq
Owing to \eqref{11.02.10},
$$\lim\limits_{\veps\rightarrow 0}\;L(G_\veps)=0,$$
and by definition of $H$ and $H_\veps$ (see (\ref{02.08.51})),
$$H_\veps\subseteq H, \ \ \  \lim\limits_{\veps\rightarrow 0}L(H\setminus H_\veps)=0,$$
\beq\label{22.1}
\int_0^1f(s)ds=\int_H f(s)ds + \int_{[0,1]\setminus H} f(s)ds
 = \int_H f(s)ds.
\eeq
Therefore, because of \eqref{11.02.11}, \eqref{22.1}, for any $\gamma\in(0,1)$
$$\begin{aligned}
&\lim\limits_{\veps\rightarrow 0} \limsup\limits_{T\rightarrow\infty} \
 \dfrac{\ln\mathbf{P}\bigl(\xi_T(\cdot)\in  \mathbb{U}_\veps(f)\bigr)}{\varphi(T)\ln\varphi(T)}  \\
&\quad \leq  -\dfrac{Q(1-\gamma)}{k}\diy\int_Hf(s)ds
+ \lim\limits_{\veps\rightarrow 0} \limsup\limits_{T\rightarrow\infty}
\dfrac{\ln\mathbf{E}\bigl[e^{B_T+N_T\ln2}
\mathbf{1}(\zeta_T\in  \mathbb{U}_\veps(f))\bigr]}{\varphi(T)\ln\varphi(T)} \\
&\quad =  -\dfrac{Q(1-\gamma)}{k}\diy\int_0^1f(s)ds
+\lim\limits_{\veps\rightarrow 0} \limsup\limits_{T\rightarrow\infty}
\dfrac{\ln\mathbf{E}\bigl[e^{B_T+N_T\ln2}\mathbf{1}(\zeta_T\in \mathbb{U}_\veps(f))\bigr]
}{\varphi(T)\ln\varphi(T)}.\end{aligned}$$
Passing to the limit $\gamma\rightarrow 0$ and using Lemma \ref{l2.4}, we get for $f\in \mathbb{V}$
\beq\label{11.02.12}
\lim\limits_{\veps\rightarrow 0} \limsup\limits_{T\rightarrow\infty} \
\dfrac{\ln\mathbf{P}\bigl(\xi_T(\cdot)\in  \mathbb{U}_\veps(f)\bigr)}{\varphi(T)\ln\varphi(T)}
\leq -\dfrac{Q}{k}\int_0^1f(s)ds-(1-l)f^+_\mathbb{D}(1). \eeq
Owing to the lower bound in \eqref{17.12.18}, and using an argument similar to the one
above, together with Lemma \ref{l2.5}, we obtain for $f\in \mathbb{V}$
\beq\label{11.02.14}
\lim\limits_{\veps\rightarrow 0} \liminf\limits_{T\rightarrow\infty} \
\dfrac{\ln\mathbf{P}\bigl(\xi_T(\cdot)\in  \mathbb{U}_\veps(f)\bigr)}{\varphi(T)\ln\varphi(T)}
\geq -\dfrac{Q}{k}\int_0^1f(s)ds-(1-l)f^+_\mathbb{D}(1). \eeq
This completes the proof of the $(\mathbb{V},\mathbb{L},\rho)$-LLDP in assertion 1).
\vskip .2cm

Now, assume condition \eqref{17.12.6}. Then bound \eqref{17.12.18}, along with Lemmas \ref{l2.4}
and \ref{l2.5}, implies that
$$\begin{aligned}
&\lim\limits_{\veps\rightarrow 0} \limsup\limits_{T\rightarrow\infty}
\dfrac{\ln\mathbf{P}\bigl(\xi_T(\cdot)\in  \mathbb{U}_\veps(f)\bigr)}{\varphi(T)\ln\varphi(T)}\\
&\qquad =\lim\limits_{\veps\rightarrow 0} \liminf\limits_{T\rightarrow\infty}
\dfrac{\ln\mathbf{P}\bigl(\xi_T(\cdot)\in  \mathbb{U}_\veps(f)\bigr)}{\varphi(T)\ln\varphi(T)}=
-(1-l)f^+_\mathbb{D}(1).\end{aligned}$$
This completes the proof of $(\mathbb{V},\mathbb{L},\rho)$-LLDP in assertion 2).

Furthermore, Lemma \ref{l4.5} implies the ET property for family $\xi_T(\,\cdot\,)$, $T>0$
and the fact that $I(f)=\infty$ for $f\in \mathbb{L}\setminus \mathbb{V}$ under any of conditions
(\ref{11.02.7}) or (\ref{17.12.6}). As a result, we get an LDP under each of conditions
\eqref{11.02.7}, \eqref{17.12.6}.
\quad $\Box$
\vspace{0.5cm}

The proofs of Lemmas \ref{l2.4} and \ref{l2.5} are based on upper and lower bounds
for the expected value
$$E=:
\mathbf{E}\bigl( e^{B_T+N_T(\zeta)\ln2}\mathbf{1}[\zeta_T\in  \mathbb{U}_\veps(f)]\bigr).$$
\vspace{0.5cm}

{\sc Proof of Lemma \ref{l2.4}.}
Given
$a\in (0,\infty )$, let $\mathbb{V}_a$ be the set of functions $f\in\mathbb{V}$ with
$0\leq f(0)\leq a$ and $\mathrm{Var}\,f\leq a$. Next, given $C\in (0,\infty )$,
define the set $\mathbb{K}_C:=\mathbb{V}_{a(C)}$ with $$a(C) :=\dfrac{3C}{1-l}.$$ According to
Lemma~\ref{l4.2}, $\mathbb{K}_C$ is compact in $(\mathbb{L},\rho)$. We write
\beq\label{17.12.19}
\begin{aligned}
&E\leq E_1+E_2,\quad\hbox{where}\\
&\quad E_1:=
\mathbf{E}\bigl(e^{B_T+N_T(\zeta)\ln2}\mathbf{1}[\zeta_T\in  \mathbb{U}_\veps(f)\cap\mathbb{K}_C]\bigr),\\
&\quad 
E_2:=
\mathbf{E}\bigl(e^{B_T+N_T(\zeta)\ln2}\,\mathbf{1}[\zeta_T\in \mathbb{K}_C^\complement]\bigr), \\
\end{aligned}\eeq
and $\mathbb{K}_C^\complement$ stands for the complement to set $\mathbb{K}_C$.

Let us upper-bound the term $E_1$. Obviously, process $\zeta(t)$ can be represented as
\beq\label{12.02.6}\zeta(t)=\zeta^+(t)-\zeta^-(t),\eeq
where $\zeta^+$ and $\zeta^-$ are independent Poisson processes of rate $1/2$, with
$$\mathbf{E} \bigl( \zeta^+(t) \bigr) =\mathbf{E}\bigl( \zeta^-(t) \bigr) =t/2.$$
Note that if $\zeta_T\in\mathbb{K}_C$ then, by virtue of \eqref{17.12.3},
for any $\gamma\in (0,1)$ and $T$ large enough, we can upper bound $B_T$ as follows
\beq\label{17.12.20}\begin{aligned}
B_T= & \sum_{i=1}^{N_T(\zeta)}\ln\bigl( \nu(\zeta(t_{i-1}),\zeta(t_i)\bigr)\\
\leq & \
 \zeta^-(T)\ln\bigl(Q\varphi(T)a(C)(1+\gamma)\bigr) +\zeta^+(T)\ln\bigl(P\varphi^l(T)a(C)
(1+\gamma)\bigr) \phantom{\sum} \\
 =: & \  B_T^-+B_T^+. \phantom{\sum}\end{aligned}\eeq

Recall that  processes $\zeta^-$ and $\zeta^+$ are independent and non-decreasing. Also note that
$N_T(\zeta) \leq a(C) \varphi(T)$. Because of this, and due to representation \eqref{17.12.20}, Lemma
\ref{c2}, and Lemma \ref{l4.6} imply that
\beq\label{18.12.1}\begin{aligned}
E_1\leq  \
e^{a(C)\varphi(T)\ln2} &\mathbf{E} \bigl( e^{B_T}\mathbf{1}[\zeta^-_T(1)\geq f^-_\mathbb{D}(1)
-\delta(\veps)]\mathbf{1}[\zeta^+_T(1) \geq f^+_\mathbb{D}(1) -\delta(\veps)] \bigr) \\
\leq  \ e^{a(C)\varphi(T)\ln2}&\mathbf{E} \bigl( e^{B_T^-}\mathbf{1}[a(C)\geq\zeta^-_T(1)\geq f^-_\mathbb{D}(1)
-\delta(\veps)] \bigr) \\
 \times &\mathbf{E} \bigl( e^{B_T^+}\mathbf{1}[\zeta^+_T(1) \geq f^+_\mathbb{D}(1) -\delta(\veps)] \bigr).
\end{aligned}\eeq
Here $\lim\limits_{\veps\rightarrow 0}\delta(\veps)=0$, and
$$
\zeta^+_T(t):=\dfrac{\zeta^+(t T)}{\varphi(T)}, \ \ \ \zeta^-_T(t):=\dfrac{\zeta^-(t T)}{\varphi(T)}.$$

Observe that, as $\rho(f,0)>0$, for $\veps>0$  small enough
\beq \label{12.10.3} f^+_\mathbb{D}(1) -\delta(\veps)>0.\eeq

By using \eqref{17.12.20}, we obtain that as
$T\rightarrow\infty$
\beq\label{02.08.53}\begin{aligned}
& \mathbf{E}\bigl(e^{B_T^-}\mathbf{1}[a(C)\geq\zeta^-_T(1)\geq f^-_\mathbb{D}(1) -\delta(\veps)]\bigr)\\
& \qquad \leq \sum\limits_{k=0}^{ \lfloor \varphi(T)a(C) \rfloor}
e^{k\ln(Q\varphi(T)a(C)(1+\gamma))}\dfrac{e^{-T/2}(T/2)^k}{k!}\\
& \qquad\leq e^{\lfloor \varphi(T)a(C) \rfloor\ln(Qa(C)(1+\gamma))}\sum\limits_{k=0}^{ \lfloor \varphi(T)a(C) \rfloor}
e^{k\ln\varphi(T)}\dfrac{e^{-T/2}(T/2)^k}{k!}\\
& \qquad=e^{o(\varphi(T)\ln\varphi(T))}\sum\limits_{k=0}^{ \lfloor \varphi(T)a(C) \rfloor}
e^{k\ln\varphi(T)}\dfrac{e^{-T/2}(T/2)^k}{k!},\end{aligned}\eeq
where $\lfloor b\rfloor$ denotes the integer part of $b$.

It follows from the Lemma \ref{l5.7} that for $T$ large enough
$$
\max\limits_{0\leq k \leq \lfloor \varphi(T)a(C) \rfloor}e^{k\ln\varphi(T)}\dfrac{e^{-T/2}(T/2)^k}{k!}=
e^{\lfloor \varphi(T)a(C) \rfloor\ln\varphi(T)}\dfrac{e^{-T/2}(T/2)^{\lfloor \varphi(T)a(C) \rfloor}}{\lfloor \varphi(T)a(C) \rfloor!}.
$$
Therefore, using (\ref{02.08.53}) and Stirling formula, we obtain that as
$T\rightarrow\infty$
\beq\label{12.02.1}\begin{aligned}
& \mathbf{E}\bigl(e^{B_T^-}\mathbf{1}[a(C)\geq\zeta^-_T(1)\geq f^-_\mathbb{D}(1) -\delta(\veps)]\bigr)\\
& \qquad \leq (\lfloor \varphi(T)a(C)\rfloor+1)
\dfrac{e^{\lfloor \varphi(T)a(C)\rfloor \ln\varphi(T)+o(\varphi(T)\ln\varphi(T))}}{\lfloor \varphi(T)a(C)\rfloor!}
=e^{o(\varphi(T)\ln\varphi(T))}.
\end{aligned}\eeq

Applying \eqref{17.12.20}, \eqref{12.10.3} and the Stirling formula, we obtain that
for $T$ large enough
\beq\label{12.02.2}\begin{aligned}
& \mathbf{E} \bigl( e^{B_T^+}\mathbf{1}[\zeta^+_T(1)\geq f^+_\mathbb{D}(1) -\delta(\veps)] \bigr) \\
& \quad\leq\sum\limits_{k=\lfloor\varphi(T)(f^+_\mathbb{D}(1) -\delta(\veps))\rfloor}^\infty
e^{k\ln(P\varphi^l(T)a(C)(1+\gamma))}\dfrac{e^{-T/2}(T/2)^k}{k!}\\
& \quad\leq\sum\limits_{k=\lfloor \varphi(T)(f^+_\mathbb{D}(1) -\delta(\veps))\rfloor}^\infty
\exp\{k\ln(P\varphi^l(T)a(C)(1+\gamma))-k\ln k+k\ln(eT/2)\}\\
& \quad\leq\sum\limits_{k=\lfloor \varphi(T)(f^+_\mathbb{D}(1) -\delta(\veps))\rfloor}^\infty
\exp\{lk\ln\varphi(T)-k\ln k+2k\ln(eT/2)\}\\
& \quad\leq\sum\limits_{k=\lfloor \varphi(T)(f^+_\mathbb{D}(1) -\delta(\veps))\rfloor}^\infty
\exp\big\{lk\ln\varphi(T)-k\ln (\lfloor\varphi(T)(f^+_\mathbb{D}(1) -\delta(\veps))\rfloor)+2k\ln(eT/2)\big\}\\
& \quad\leq\sum\limits_{k=\lfloor\varphi(T)(f^+_\mathbb{D}(1) -\delta(\veps))\rfloor}^\infty
\exp\{-(1-l)k\ln \varphi(T)+3 k \ln (eT/2)\}\\
& \quad=\sum\limits_{k=\lfloor\varphi(T)(f^+_\mathbb{D}(1) -\delta(\veps))\rfloor}^\infty
\exp\big\{-k\big((1-l)\ln\varphi(T)-3  \ln (eT/2)\big)\big\}\\
& \quad=\dfrac{\exp\big\{-\lfloor\varphi(T)(f^+_\mathbb{D}(1) -\delta(\veps))\rfloor\big((1-l)\ln\varphi(T)-3
\ln (eT/2)\big)\big\}}{1-\exp\big\{-(1-l)\ln\varphi(T)+3  \ln (eT/2)\big\}}.\end{aligned}\eeq

From bounds \eqref{18.12.1}, \eqref{12.02.1}, \eqref{12.02.2} we get that for $T$ large enough
\beq\label{12.02.3}
E_1\leq\exp\big\{ -(1-l)(f^+_\mathbb{D}(1) -\delta(\veps))\varphi(T)\ln \varphi(T)
+o(\varphi(T)\ln\varphi(T))\big\} \eeq
Further, Lemma \ref{l4.5} implies that for $T$ large enough
\beq\label{12.02.4}
E_2\leq\exp\big\{ -C\varphi(T)\ln\varphi(T)+o(\varphi(T)\ln\varphi(T))\big\}.
\eeq
Choosing $C>(1-l)(f^+_\mathbb{D}(1) -\delta(\veps))$ and using inequalities
\eqref{17.12.19}, \eqref{12.02.3}, \eqref{12.02.4}, we obtain
$$\begin{aligned}
&\lim\limits_{\veps\rightarrow 0} \limsup\limits_{T\rightarrow\infty}\;
\dfrac{\ln E}{\varphi(T)\ln\varphi(T)} \\
&\quad\leq\lim\limits_{\veps\rightarrow 0} \limsup\limits_{T\rightarrow\infty}\;
\dfrac{\ln \left( 2e^{-(1-l)(f^+_\mathbb{D}(1) -\delta(\veps))\varphi(T)\ln \varphi(T)+o(\varphi(T)\ln\varphi(T))}\right)
 }{\varphi(T)\ln\varphi(T)} \\
&\quad\leq -\lim\limits_{\veps\rightarrow 0}(1-l)(f^+_\mathbb{D}(1) -\delta(\veps))=-(1-l)f^+_\mathbb{D}(1).
\quad\Box \phantom{\dfrac{1}{\varphi(T)\ln\varphi(T)}}\end{aligned}$$
\vspace{0.5cm}

\noindent
{\sc Proof of Lemma \ref{l2.5}.} Here the goal is to establish a lower bound for $E$.
As usually, obtaining lower bounds is a more difficult task. Let us outline the
idea of the proof. The main step is to extract from the event $\big\{\zeta_T\in
\mathbb{U}_\veps(f)\big\}$ a smaller event:
$$\{ \zeta_T\in  \mathbb{U}_\veps(f) \} \supset \left\{\int_0^\Delta|\zeta_T(t)- \widetilde{g}_\veps(t)|dt
<\veps/4, \  \sup_{t\in [\Delta,1]}|\zeta_T(t)- \widetilde{g}_\veps(t)| < \veps/8 \right\}.$$
Here $\Delta$ is a constant that depends on $\veps$, and $\widetilde{g}_\veps$ is a continuous function
such that
\begin{enumerate}
    \item $\widetilde{g}_\veps$ is close to $f$ in $\rho$ metric (in the proof $\rho(f, \widetilde{g}_\veps)
< \frac{3\veps}{4}$).
    \item The variation $\mathrm{Var}\,\widetilde{g}_\veps$ is close to the variation of $f$.
    \item $\widetilde{g}_\veps$ is equal to a small constant $\delta$ on the interval $[0,\Delta]$ (in
the proof $\delta = \veps/4$) and $\widetilde{g}_\veps$ is greater than $\delta$ on $[\Delta, 1]$.
\end{enumerate}
Then the expected value $E$ will be lower-bounded by a product $E_+E_-$ (see
(\ref{13.10.2})), where $E_+$ (respectively, $E_-$) controls the variations of
$\zeta^+$ and $\widetilde{g}_\veps^+$ (respectively, $\zeta^-$ and $\widetilde{g}_\veps^-$). Finally, the
quantity $E_+$ will give us the bound claimed in the lemma \eqref{13.10.5}).

\vspace{0.5cm}
Let us proceed with the formal proof. First, consider the case where
$\mathrm{Var}\,f_{\mathbb{D}}>0$.  We start by proving
the existence of function $\widetilde{g}_\veps$. We construct $\widetilde{g}_\veps$  by
using an auxiliary function $g$ (see below). From the point of view of future arguments, it is convenient
to set:
$$2\Delta:=\sup\Bigl\{t: \mathrm{Var}\,f_{\mathbb{D}}\leq\dfrac{\veps}{2}\Bigr\}.$$
Observe that since $\mathrm{Var}\,f_{\mathbb{D}}>0$, we have $2\Delta<1$ for $\veps>0$
small enough and, as $f_{\mathbb{D}}$ is right-continuous at $0$, we also have that
$2\Delta>0$.

Define
$$g(t):=\left\{ \begin{array}{ll} \dfrac{\veps}{4}\vee f(0), \ t\in [0,2\Delta),\\
f(t)+\dfrac{\veps}{2},   \  t\in [2\Delta,1],
\end{array} \right.$$
Recall that $a\vee b=\max(a,b)$. It is easy to see that $\rho(f,g)\leq\dfrac{\veps}{2}$.
Note that function $g$ is convenient because it does not vanish on $[0,1]$.

Let us decompose $g_\mathbb{D}$ into an increasing and a decreasing component:
$$\begin{aligned}
g_\mathbb{D}(t) &= g^+_\mathbb{D}(t)-g^-_\mathbb{D}(t), \ \
\mathrm{Var}\,g_{\mathbb{D}}=\mathrm{Var}\,g^+_{\mathbb{D}}+\mathrm{Var}\,g^-_{\mathbb{D}} \\
g^+_\mathbb{D}(t) &= \frac{\veps}{4}\vee f(0), \ \text{for}\ t\in [0,2\Delta),\\
 g^-_\mathbb{D}(t) &=0,  \ \text{for} \
t\in [0,2\Delta).\end{aligned}
$$
From the definition of the constant $2\Delta$ it follows that
\beq \label{13.10.8}
0\leq\mathrm{Var}\,f^+_{\mathbb{D}}-\mathrm{Var}\,g^+_{\mathbb{D}}\leq \dfrac{\veps}{2}.
\eeq
Functions $g^+_\mathbb{D}$, $g^-_\mathbb{D}$ are monotone and continuous
on $[0,2\Delta)$ and left-continuous at the end point $1$. Also,
$$\inf\limits_{t\in[0,2\Delta)}(g^+_\mathbb{D}(t)-g^-_\mathbb{D}(t))=\dfrac{\veps}{4}\vee f(0), \
\inf\limits_{t\in[2\Delta,1]}(g^+_\mathbb{D}(t)-g^-_\mathbb{D}(t))\geq\dfrac{\veps}{2}.$$

Hence, there exist monotone continuous functions
$\widetilde{g}^+_\veps$, $\widetilde{g}^-_\veps$ such that
$$\begin{aligned}
& \widetilde{g}^+_\veps(t) = g^+_\mathbb{D}(t)=\dfrac{\veps}{4}\vee f(0), \
\widetilde{g}^-_\veps(t)=g^-_\mathbb{D}(t)=0, \ \text{if} \ t\in [0,\Delta],\\
& \widetilde{g}^+_\veps(t) \geq g^+_\mathbb{D}(t), \
\widetilde{g}^-_\veps(t)\leq g^-_\mathbb{D}(t), \  \text{if}  \
t\in (\Delta,1), \phantom{\frac{\veps}{4}} \\
&\widetilde{g}^+_\veps(1)  =g^+_\mathbb{D}(1), \
\widetilde{g}_\veps^-(1)=g^-_\mathbb{D}(1), \text{ and}\phantom{\frac{\veps}{4}} \\
& \rho(\widetilde{g}^+_\veps,g^+_\mathbb{D})+\rho(\widetilde{g}^-_\veps,g^-_\mathbb{D})
<\frac{\veps}{4}.\end{aligned}$$
Let $\widetilde{g}_\veps:=\widetilde{g}^+_\veps-\widetilde{g}^-_\veps$. Then
$\rho(\widetilde{g}_\veps,g)<\dfrac{\veps}{4}$ and
$$\inf\limits_{t\in[0,1]}\widetilde{g}_\veps(t)=\inf\limits_{t\in[0,1]}\bigl(\widetilde{g}^+_\veps(t)
-\widetilde{g}^-_\veps(t)\bigr)
\geq\inf\limits_{t\in[0,1]}\bigl(g^+_\mathbb{D}(t)-g^-_\mathbb{D}(t)\bigr)\geq\frac{\veps}{4}.$$

Using decomposition \eqref{12.02.6}, we have for $T$ large enough
\beq\label{12.02.5}
\begin{aligned}
&E:= \ \mathbf{E}\left(e^{B_T+N_T(\zeta)\ln2}
\mathbf{1}[\zeta_T\in  \mathbb{U}_\veps(f)]\right) \phantom{\dfrac{\veps}{8}}\\
&\geq \ \mathbf{E}\Bigl(e^{B_T}
\mathbf{1}\Big[\sup\limits_{t\in[0,\Delta)}\zeta_T(t)<\widetilde{g}_\veps(\Delta),
\sup\limits_{t\in[\Delta,1]}|\zeta_T(t)-\widetilde{g}_\veps(t)|<\dfrac{\veps}{8}\Big]\Bigr)\\
&\geq\ \mathbf{E}\Bigl(e^{B_T}\mathbf{1}\Big[\zeta^-(\Delta T)=0,  \zeta^+(\Delta T)
=\lfloor \varphi(T)\widetilde{g}_\veps(\Delta)\rfloor,
\sup\limits_{t\in[\Delta,1]}|\zeta_T(t)-\widetilde{g}_\veps(t)|<\dfrac{\veps}{8}\Big]\Bigr).
\end{aligned}\eeq

Denote
$$\begin{aligned}
\mathbb{W}_1 &:=\Bigl\{
\sup\limits_{t\in[\Delta,1]} \left| \zeta^-_T(t)-\widetilde{g}^-_\veps(t) \right|<\dfrac{\veps}{16}\Bigr\},\\
\mathbb{S}_1 &:=\left\{\zeta^-(\Delta T)=0\right\}, \phantom{\dfrac{\veps}{8}} \\
\mathbb{W}_2 &:=\Bigl\{
\sup\limits_{t\in[\Delta,1]}\left|\zeta^+_T(t)-\zeta_T^+(\Delta)-(\widetilde{g}^+_\veps(t)
-\widetilde{g}^+_\veps(\Delta))\right| <\dfrac{\veps}{16}\Bigr\},\\
\mathbb{S}_2&:=\left\{\zeta^+(\Delta T)=\lfloor \varphi(T)\widetilde{g}_\veps(\Delta)\rfloor \right\}.
\end{aligned}$$
From bound \eqref{12.02.5} it follows that
\beq\label{13.10.1}
E\geq\mathbf{E}\Big[e^{B_T}\mathbf{1}\big( \mathbb{W}_1\cap \mathbb{W}_2\cap \mathbb{S}_1\cap \mathbb{S}_2\big)\Big].\eeq

Let us first estimate $B_T$ from below. According to the definition of $\nu$
in the lemma~\label{l2.3} the sum $B_T$ can be always separated into two sums: the one
over negative jumps, $\sum_{(-)}$ and the one over positive jumps, $\sum_{(+)}$:
$$\begin{array}{cl}
B_T &=  \sum\limits_{i=1}^{N_T(\zeta)}\ln(\nu(\zeta(t_{i-1}),\zeta(t_i)) \\
&= \sum\limits_{i=1}^{N_T(\zeta)}\ln(\nu(\zeta(t_{i-1}),\zeta(t_i))
\mathbf{1}[\zeta(t_{i-1})>\zeta(t_i)]\\
&\qquad +\sum\limits_{i=1}^{N_T(\zeta)}\ln(\nu(\zeta(t_{i-1}),\zeta(t_i))
\mathbf{1}[\zeta(t_{i-1}) <\zeta(t_i)] =:  \sum_{(-)} + \sum_{(+)}.\end{array}$$

The lower bound for $\sum_{(-)}$ is constructed in the following way. Note that, according to
$\mathbb{S}_1$, there are no negative jumps of the process $\zeta$ during the time interval
$[0,T\Delta]$ and, according to
$\mathbb{W}_1\cap\mathbb{W}_2$, the process belongs to an $\veps/8$ neighborhood (in the
uniform metric) of the function  $\widetilde{g}_\veps$. Then, due to (\ref{17.12.3}),
for any $\gamma\in(0,1)$ and $T$ sufficiently large we obtain
\beq\label{-}\begin{array}{l}
\sum_{(-)} \geq\zeta^-(T)\ln\Bigl(\frac{\veps Q\varphi(T)(1-\gamma)}{8}\Bigr).
\end{array}\eeq

Let us bound sum $\sum_{(+)}$ from below. Note that for each $r>1$ any  trajectory from the event
$\mathbb{S}_1\cap \mathbb{S}_2\cap \mathbb{W}_2$ has
$\zeta^+(T)-\lfloor\frac{\varphi(T)\widetilde{g}_\veps(\Delta)}{r}\rfloor$
positive jumps, when when the trajectory is not lower than  $\veps/8r$.
Thus, on these jumps for any $\gamma\in(0,1)$ and $T$ sufficiently large we have
$$
\nu(\zeta(t_{i-1}),\zeta(t_i))\geq \frac{\veps P\varphi^l(T)(1-\gamma)}{8r}.
$$

On the remaining $\lfloor\frac{\varphi(T)\widetilde{g}_\veps(\Delta)}{r}\rfloor$ positive jumps,
$\nu(\zeta(t_{i-1}),\zeta(t_i))\geq \lambda_{\min} 
$; here $\lambda_{\min}:=\min\limits_{x\in \mathbb{Z}^+}\lambda(x)$.

Finally,
\beq\label{+}
\sum_{(+)}
\geq\Bigl(\zeta^+(T)-\Bigl\lfloor \dfrac{\varphi(T)\widetilde{g}_\veps(\Delta)}{r}\Bigr\rfloor \Bigr)
\ln\Bigl(\dfrac{\veps P\varphi^l(T)(1-\gamma)}{8r}\Bigr)
+\Bigl\lfloor \dfrac{\varphi(T)\widetilde{g}_\veps(\Delta)}{r} \Bigr\rfloor\ln \lambda_{\min}.
\eeq
The introduced new parameter $r$ will further tend to infinity, $r\rightarrow\infty$.

Thus, on the event $\mathbb{W}_1\cap\mathbb{W}_2\cap \mathbb{S}_1\cap \mathbb{S}_2$ due to \eqref{-} and
\eqref{+} we obtain
that for any $\gamma\in(0,1)$, $r>1$ and $T$ large enough,
\beq\label{12.02.8}
\begin{aligned}
B_T= & \sum\limits_{i=1}^{N_T(\zeta)}\ln(\nu(\zeta(t_{i-1}),\zeta(t_i))) \\
\geq & \
 \zeta^-(T)\ln\Bigl(\dfrac{\veps Q\varphi(T)(1-\gamma)}{8}\Bigr) + \Bigl\lfloor \dfrac{\varphi(T)\widetilde{g}_\veps(\Delta)}{r} \Bigr\rfloor\ln \lambda_{\min} \\
 & +\Bigl(\zeta^+(T)-\Bigl\lfloor \dfrac{\varphi(T)\widetilde{g}_\veps(\Delta)}{r}\Bigr\rfloor \Bigr)
\ln\Bigl(\dfrac{\veps P\varphi^l(T)(1-\gamma)}{8r}\Bigr) \\
\geq & \  \zeta^-(T)\ln(J\varphi(T)) + \Bigl\lfloor \dfrac{\varphi(T)\widetilde{g}_\veps(\Delta)}{r}\Bigr\rfloor \ln \lambda_{\min}\\
&+\Bigl(\zeta^+(T)-\Bigl\lfloor \dfrac{\varphi(T)\widetilde{g}_\veps(\Delta)}{r}\Bigr\rfloor \Bigr)\ln(J\varphi^l(T)), \end{aligned}\eeq
where
$J:=\min\Bigl(\dfrac{\veps P(1-\gamma)}{8r},\dfrac{\veps Q(1-\gamma)}{8}\Bigr)$. \vskip .3cm

Owing to \eqref{13.10.1}, \eqref{12.02.8} and independence of processes $\zeta^+$ and $\zeta^-$, we get
\beq\label{13.10.2}
\begin{aligned}
E\geq & \mathbf{E}\Bigl(\exp\Bigl\{ \big(\zeta^+(T)-\Bigl\lfloor \frac{\varphi(T)
\widetilde{g}_\veps(\Delta)}{r}\Bigr\rfloor \big)\ln(J\varphi^l(T)) \\
& \qquad \qquad +\Bigl\lfloor \frac{\varphi(T)\widetilde{g}_\veps(\Delta)}{r}\Bigr\rfloor \ln \lambda_{\min}\Bigr\}
\mathbf{1}\big(\mathbb{S}_2\cap \mathbb{W}_2\big)\Bigr)\\
& \times\mathbf{E}\Big(\exp\big\{\zeta^-(T)\ln(J\varphi(T))\big\}\mathbf{1}\big(\mathbb{S}_1\cap \mathbb{W}_1\big)\Big)
=:E_+E_-. \phantom{\frac{\varphi(T)\widetilde{g}_\veps(\Delta)}{r}}
\end{aligned}\eeq

Let us lower-bound the value $E_+$. Consider a partition
$\Delta=t_0<t_1<\dots<t_m=1$ such that
$$
\max\limits_{i=1,\dots,m}\bigl(\widetilde{g}^+_\veps(t_i)-\widetilde{g}^+_\veps(t_{i-1})\bigr)\leq\dfrac{\veps}{32}
\ \ \ \text{and} \ \ \ \min\limits_{i=1,\dots,m}\bigl(\widetilde{g}^+_\veps(t_i)-\widetilde{g}^+_\veps(t_{i-1})\bigr)>0.$$

By virtue of independence of increments in process $\zeta^+$, for $T$ large enough
\beq\label{13.10.4}
\begin{aligned}
E_+\geq \mathbf{E} \Bigl( & \exp\Bigl\{ \bigl(\zeta^+(T)-\Bigl\lfloor \frac{\varphi(T)\widetilde{g}_\veps(\Delta)}{r}\Bigr\rfloor \bigr)
\ln (J\varphi^l(T)) + \Bigl\lfloor \frac{\varphi(T)\widetilde{g}_\veps(\Delta)}{r}\Bigr\rfloor \ln (\lambda_{\min}) \Bigr\}
\\
& \times\mathbf{1}\bigl[ \zeta^+(\Delta T) = \lfloor \varphi(T)\widetilde{g}_\veps(\Delta) \rfloor \bigr]  \\
& \times\prod\limits_{i=1}^m\mathbf{1} \bigl[\zeta^+(Tt_i)-\zeta^+(Tt_{i-1})
=\lfloor (\widetilde{g}^+_\veps(t_i)-\widetilde{g}^+_\veps(t_{i-1}))\varphi(T)\rfloor \bigr] \Bigr)\\
\geq \mathbf{E}\Bigl( & \exp\Bigr\{ \lfloor\varphi(T)\widetilde{g}^+_\veps(1)\rfloor \ln(J\varphi^l(T))-
2\Bigl\lfloor \frac{\varphi(T)\widetilde{g}_\veps(\Delta)}{r}\Bigr\rfloor \ln(J\varphi^l(T))\Bigr\} \\
& \times \mathbf{1} \left[\zeta^+(\Delta T)= \lfloor \varphi(T)\widetilde{g}_\veps(\Delta) \rfloor \right] \\
&\times\prod\limits_{i=1}^m\mathbf{1} \left[ \zeta^+(Tt_i)-\zeta^+(Tt_{i-1})
= \lfloor (\widetilde{g}^+_\veps(t_i)-\widetilde{g}^+_\veps(t_{i-1}))\varphi(T)\rfloor \right] \Bigr) \\
= &\  \exp\Bigr\{ \lfloor\varphi(T)\widetilde{g}^+_\veps(1)\rfloor \ln(J\varphi^l(T))-
2\Bigl\lfloor \frac{\varphi(T)\widetilde{g}_\veps(\Delta)}{r}\Bigr\rfloor \ln(J\varphi^l(T))\Bigr\}
\\
& \times
\mathbf{P}\left(\zeta^+(\Delta T)=\lfloor \varphi(T)\widetilde{g}_\veps(\Delta)\rfloor \right)\\
&\times\prod\limits_{i=1}^m\mathbf{P}\left(\zeta^+(Tt_i)-\zeta^+(Tt_{i-1})
=\lfloor (\widetilde{g}^+_\veps(t_i)-\widetilde{g}^+_\veps(t_{i-1}))\varphi(T)\rfloor \right).
\end{aligned}
\eeq
With the help of the Stirling formula, we get that for $T$ large enough
\beq \label{12.02.9}\begin{aligned}
 \mathbf{P}&\left(\zeta^+(\Delta T)=\lfloor \varphi(T)\widetilde{g}_\veps(\Delta)\rfloor \right) \\
 & \times \prod\limits_{i=1}^m\mathbf{P} \left( \zeta^+(Tt_i)-\zeta^+(Tt_{i-1})
= \lfloor (\widetilde{g}^+_\veps(t_i)-\widetilde{g}^+_\veps(t_{i-1}))\varphi(T) \rfloor \right)\\
= & \dfrac{e^{-T\Delta/2}(T\Delta/2)^{\lfloor \widetilde{g}_\veps(\Delta)\varphi(T) \rfloor }}
{\lfloor \widetilde{g}_\veps(\Delta)\varphi(T)\rfloor !}
\prod\limits_{i=1}^m \dfrac{e^{-T(t_i-t_{i-1})/2}(T(t_i-t_{i-1})/2)^{ \lfloor (\widetilde{g}^+_\veps(t_i)-\widetilde{g}^+_\veps(t_{i-1}))\varphi(T)\rfloor }}
{\lfloor (\widetilde{g}^+_\veps(t_i)-\widetilde{g}^+_\veps(t_{i-1}))\varphi(T)\rfloor !}\\
\geq &
\prod\limits_{i=1}^m \exp\left\{-\dfrac{T(t_i-t_{i-1})}{2}-
(\widetilde{g}^+_\veps(t_i)-\widetilde{g}^+_\veps(t_{i-1}))\varphi(T)\ln((\widetilde{g}^+_\veps(t_i)-\widetilde{g}^+_\veps(t_{i-1}))\varphi(T))\right\}\\
&\times \exp\left\{ -\frac{T\Delta}{2}-\widetilde{g}_\veps(\Delta)\varphi(T)\ln(\widetilde{g}_\veps(\Delta)
\varphi(T)) \right\}   \phantom{ -\dfrac{T(t_i-t_{i-1})}{2} } \\
\geq &
\prod\limits_{i=1}^m \exp\left\{-\dfrac{T(t_i-t_{i-1})}{2}-
(\widetilde{g}^+_\veps(t_i)-\widetilde{g}^+_\veps(t_{i-1}))\varphi(T)\ln(\widetilde{g}^+_\veps(1)\varphi(T))\right\}\\
&\times \exp\left\{ -\frac{T\Delta}{2}-\widetilde{g}^+_\veps(\Delta)\varphi(T)\ln(\widetilde{g}^+_\veps(1)\varphi(T))\right\} \\
\geq & \exp\left\{-T-
\widetilde{g}^+_\veps(1)\varphi(T)\ln(\widetilde{g}^+_\veps(1)\varphi(T))\right\}. \phantom{ -\dfrac{T(t_i-t_{i-1})}{2} }
\end{aligned}
\eeq
From bounds \eqref{13.10.4}, \eqref{12.02.9} it follows that for $T$ large enough
\beq \label{13.10.5}\begin{aligned}
& E_+\geq
\exp\Bigl\{-\widetilde{g}^+_\veps(1)\varphi(T)(1-l)\ln \varphi(T)\\
& \qquad\qquad - T- \varphi(T)\widetilde{g}^+_\veps(1)|
\ln (J)|
 -3\Bigl\lfloor \dfrac{\varphi(T)\widetilde{g}_\veps(\Delta)}{r}\Bigr\rfloor \ln(J\varphi^l(T)) \Bigr\}.
\end{aligned}\eeq

Next, let us lower-bound the quantity $E_-$. If $\widetilde{g}^-_\veps(1)=0$ then
\beq \label{13.10.6}
E_-\geq \mathbf{P}(\zeta^-(T)=0)=e^{-T/2}.
\eeq
If $\widetilde{g}^-_\veps(1)>0$ then we consider the partition
$\Delta=t_0<t_1<\dots<t_m=1$ such that
$$
\max\limits_{i=1,\dots,m}\bigl(\widetilde{g}^-_\veps(t_i)-\widetilde{g}^-_\veps(t_{i-1})\bigr)
\leq\dfrac{\veps}{32}
\ \ \ \text{and} \ \ \ \min\limits_{i=1,\dots,m}\bigl(\widetilde{g}^-_\veps(t_i)-\widetilde{g}^-_\veps(t_{i-1})\bigr)>0.
$$
By using the independence of increments in process $\zeta^-$ and the Stirling formula,
for $T$ large enough
\beq\label{13.10.7}\begin{aligned}
E_-\geq &\  \mathbf{E} \Bigl( e^{\zeta^-(T)\ln(J\varphi(T))}\mathbf{1}(\zeta^-(\Delta T)=0)\\
& \times\prod\limits_{i=1}^m\mathbf{1} \left[ \zeta^-(Tt_i)-\zeta^-(Tt_{i-1})
=\lfloor (\widetilde{g}^-_\veps(t_i)-\widetilde{g}^-_\veps(t_{i-1}))\varphi(T) \rfloor \right] \Bigr) \\
\geq & \ \exp\left\{ \varphi(T)\widetilde{g}_\veps^-(1)\ln(J\varphi(T))-
m\ln(J\varphi(T))\right\} \mathbf{P}(\zeta^-(\Delta T)=0)\\
&\times\prod\limits_{i=1}^m\mathbf{P}\left( \zeta^-(Tt_i)-\zeta^-(Tt_{i-1})
=\lfloor (\widetilde{g}^-_\veps(t_i)-\widetilde{g}^-_\veps(t_{i-1}))\varphi(T)\rfloor \right)\\
\geq &\  \exp\left\{ \varphi(T)\widetilde{g}_\veps^-(1)\ln(J\varphi(T))-
m\ln(J\varphi(T))\right\} \\
& \times \exp\left\{-T-
\widetilde{g}^-_\veps(1)\varphi(T)\ln(\widetilde{g}^-_\veps(1)\varphi(T))\right\} \phantom{\prod\limits_{i=1}^m}\\
=& \ \exp\left\{-T-\widetilde{g}^-_\veps(1)\varphi(T)\ln \widetilde{g}^-_\veps(1)+
\widetilde{g}^-_\veps(1)\varphi(T)\ln J-
m\ln(J\varphi(T))\right\}. 
\end{aligned}\eeq

Owing to \eqref{13.10.1} and \eqref{13.10.5}\,--\,\eqref{13.10.7} we obtain that for any $r>1$
and $\veps>0$ small enough
$$\liminf\limits_{T\rightarrow\infty}\;\dfrac{\ln E}{\varphi(T)\ln\varphi(T)}
\geq
-\widetilde{g}^+_\veps(1)(1-l) - \dfrac{3l\widetilde{g}_\veps(\Delta)}{r}.
$$
Passing to the limit $r\rightarrow\infty$ yields
$$\liminf\limits_{T\rightarrow\infty}\;\dfrac{\ln E}{\varphi(T)\ln\varphi(T)}
\geq -\widetilde{g}^+_\veps(1)(1-l).$$
By definition, $\widetilde{g}^+_\veps(1)=g^+_\mathbb{D}(1)$. Also, by virtue of \eqref{13.10.8},
$|g^+_\mathbb{D}(1)-f^+_\mathbb{D}(1)|\leq\frac{\veps}{2}$. This gives
$$\lim\limits_{\veps\rightarrow 0}
\liminf\limits_{T\rightarrow\infty}\dfrac{\ln E}{\varphi(T)\ln\varphi(T)}
\geq -(1-l)f^+_\mathbb{D}(1).
$$
Therefore, Lemma \ref{l2.5} has been proven when $\mathrm{Var}\,f_{\mathbb{D}}>0$.

In the case where $\mathrm{Var}\,f_{\mathbb{D}}=0$, we have $f_{\mathbb{D}}(t)=f(0)>0$,
$t\in[0,1]$. It is easy to see that $\bigl\{\zeta_T\in  \mathbb{U}_\veps(f) \bigr\}\supseteq\mathbb{O}$
\ where
$$\begin{aligned}
&\mathbb{O}:= \Bigl\{ \zeta^+\Bigl(\frac{\veps}{2f(0)}T\Bigr)=\lfloor f(0)\varphi(T)\rfloor, \\
&\qquad\qquad\qquad\zeta^+(T)-\zeta^+\Bigl(\frac{\veps}{2f(0)}T\Bigr)=0,\zeta^-(T)=0\Bigr\}.\end{aligned}$$
The rest of the proof is reduced to a lower bound for  $B_T$ on event $\mathbb{O}$ which essentially repeats
the above argument. For brevity, we omit it from the paper.  \quad $\Box$

\section{Discussion}\label{Discussion}

It is instructive to discuss Theorems \ref{th2.1} and \ref{th2.2}  in connection with a question mentioned
in Section \ref{Sect1}: "why under
condition \eqref{17.12.5} we get only an LLDP whereas \eqref{11.02.7}  or \eqref{17.12.6} lead to an LDP?".
Consider an example where
$$\lambda(x)=P, \ \ \ \mu(x)=Q x. $$

In this case one can write down a probability distribution for process $\xi$ at the time point $T$ explicitly \cite{Gned}
\beq\label{14.10.2}
\mathbf{P}(\xi(T)=x)=\dfrac{(a(T))^x}{x!}e^{-a(T)}, \ \ \ x\in \mathbb{Z}^+,\eeq
where
$$a(T)=\dfrac{P}{Q}(1-e^{-QT}).$$
As follows from \eqref{14.10.2}, if $f(1)>0$ then
$$\lim\limits_{\veps\rightarrow 0}\lim\limits_{T\rightarrow\infty}\dfrac{1}{\varphi(T)\ln\varphi(T)}
\ln\mathbf{P}\bigl( \xi_T(1)\in [f(1)-\veps,f(1)+\veps] \bigr) =-f(1)$$
under any of conditions \eqref{17.12.5}--\eqref{17.12.6}.

Consequently, under condition \eqref{17.12.5} the normalizing function $\psi(T)=\varphi(T)\ln\varphi(T)$
in the LDP on the state space $\mathbb{Z}^+$ is different from function $\psi(T)=T\varphi(T)$
figuring in the LLDP on the functional
space $\mathbb{L}$. In other words, for $f\not\equiv 0$ we have that
$$\lim\limits_{\veps\rightarrow 0}\lim\limits_{T\rightarrow\infty}\dfrac{1}{\varphi(T)\ln\varphi(T)}\ln
\dfrac{\mathbf{P}\bigl(\xi_T\in  \mathbb{U}_\veps(f)\bigr)}{\mathbf{P}\bigl( \xi_T(1)\in [f(1)-\veps,f(1)
+\veps]\bigr)}=-\infty.$$
This is why under condition \eqref{17.12.5} the family of processes $\xi_T(\,\cdot\,)$, $T>0$,
lacks the ET property in any reasonable functional space.

However, under condition \eqref{11.02.7} or \eqref{17.12.6} the normalizing functions coincide, and we
manage to get an LDP in the functional space $(\mathbb{L},\rho)$ as stated in Theorem \ref{th2.2}.

Also, note that Theorem \ref{th2.2} allows us to get a rough asymptotic for the probability that a trajectory
of $\xi_T$ crosses a level $a>0$. Indeed, with the help of \eqref{17.12.8} and an argument similar to the
one used in the proof of \eqref{12.02.2} we have that under any of conditions \eqref{11.02.7} and \eqref{17.12.6}
$$
\begin{aligned}
\limsup\limits_{T\rightarrow\infty} & \frac{1}{\varphi(T)\ln\varphi(T)}\ln\mathbf{P}\Bigl(\sup\limits_{t\in[0,1]}
\xi_T(t)\geq a\Bigr)\\
&\leq\limsup\limits_{T\rightarrow\infty} \frac{1}{\varphi(T)\ln\varphi(T)}\ln\Bigl(
e^T \mathbf{E}\Bigl(  e^{-A_T(\zeta)}e^{B_T(\zeta)+N_T(\zeta)\ln2}
\mathbf{1}\big( \sup\limits_{t\in[0,1]}\zeta_T(t)\geq a \big)
\Bigr)\Bigr)\\
& \leq\limsup\limits_{T\rightarrow\infty}\frac{1}{\varphi(T)\ln\varphi(T)}\ln
\mathbf{E}\Bigl(e^{B_T(\zeta)+N_T(\zeta)\ln2} \mathbf{1}\bigl( \zeta^+_T(1)\geq a \bigr) \Bigr)\leq-(1-l)a.\end{aligned}
$$
Since process $\xi_T$ is c\`adl\`ag, and the set
$$\Bigl\{f\in \mathbb{L}:\text{ess}\sup\limits_{t\in[0,1]}f(t)> a\Bigr\}$$
is open, Theorem \ref{th2.2} implies that under \eqref{11.02.7} or \eqref{17.12.6}
$$\begin{aligned}
\liminf\limits_{T\rightarrow\infty} & \dfrac{1}{\varphi(T)\ln\varphi(T)}
\ln\mathbf{P}\Bigl(\sup\limits_{t\in[0,1]}\xi_T(t)\geq a\Bigr) \\
&\geq\liminf\limits_{T\rightarrow\infty}\frac{1}{\varphi(T)\ln\varphi(T)}
\ln\mathbf{P}\Bigl(\sup\limits_{t\in[0,1]}\xi_T(t)> a\Bigr)\\
& = \liminf\limits_{T\rightarrow\infty}\dfrac{1}{\varphi(T)\ln\varphi(T)}\ln\mathbf{P}
\Bigl(\text{ess}\sup\limits_{t\in[0,1]}\xi_T(t)> a\Bigr)\\
& \geq -\inf\limits_{f:\,\text{ess}\sup\limits_{t\in[0,1]}f(t)> a}I(f)=
-\inf\limits_{f:\,\text{ess}\sup\limits_{t\in[0,1]}f(t)> a}(1-l)f_\mathbb{D}^+(1)=-(1-l)a.
\end{aligned}$$
Thus,
$$\lim\limits_{T\rightarrow\infty}\dfrac{1}{\varphi(T)\ln\varphi(T)}
\ln\mathbf{P}\Bigl(\sup\limits_{t\in[0,1]}\xi_T(t)\geq a\Bigr)=-(1-l)a.$$

\section  {Appendix}

In this section, we establish some auxiliary assertions. \vskip .5cm

\begin{lemma} \label{l4.2} For any fixed $C>0$ the set $\mathbb{V}_C$ is a compact in
$(\mathbb{L},\rho)$.
\end{lemma} \vskip .5cm

\noindent
{\sc Proof.} The Helly theorem \cite{Nat}
implies that from every sequence $f_n\in\mathbb{V}_C$ one can extract a sub-sequence
$f_{n_k}$ convergent as $k\rightarrow\infty$ almost surely to some $f\in\mathbb{V}_C$.
Applying the Lebesgue dominated convergence theorem yields
$$\lim\limits_{k\rightarrow\infty}\int\limits_0^1|f_{n_k}(t)-f(t)|dt=
\int\limits_0^1\lim\limits_{k\rightarrow\infty}|f_{n_k}(t)-f(t)|dt=0.\quad\Box$$
\vspace{0.5cm}

Let $\mathbb{M}=\mathbb{M}[0,1]$ denote the set of non-decreasing functions on $[0,1]$.

\begin{lemma} \label{l4.3} Suppose that function $f\in \mathbb{V}$ is represented as
$$f(t)=g_1(t)-g_2(t),$$
where $g_1, g_2\in \mathbb{M}$. Then for any $0\leq t_1 < t_2 \leq 1$
$$g_1(t_2)-g_1(t_1)\geq f^+(t_2)-f^+(t_1).$$
\end{lemma}

\noindent
{\sc Proof.} Assume the opposite, then there exist
$0\leq t_1 < t_2 \leq 1$ such that $g_1(t_2)-g_1(t_1)<f^+(t_2)-f^+(t_1)$. Observe that in
that case we will also have $g_2(t_2)-g_2(t_1)<f^-(t_2)-f^-(t_1)$.

Let $\text{Var}_{[t_1,t_2]}$ stand for the variation on $[t_1.t_2]$.
Since the variation of a sum does not exceed the sum of the variations, we obtain
$$\text{Var}_{[t_1,t_2]}g_1+\text{Var}_{[t_1,t_2]}g_2=
g_1(t_2)-g_1(t_1)+g_2(t_2)-g_2(t_1)\geq\text{Var}_{[t_1,t_2]}f.$$
On the other hand,
$$
g_1(t_2)-g_1(t_1)+g_2(t_2)-g_2(t_1)<f^+(t_2)-f^+(t_1)+f^-(t_2)-f^-(t_1)=\text{Var}_{[t_1,t_2]}f.
$$
The contradiction completes the proof.\qquad $\Box$
\vspace{0.5cm}

Let $\mathbb{K}$ be a compact set in $(\mathbb{L},\rho )$.
Consider a family of functions $u_T(t)$, $t\in[0,1]$, $T>0$, such that
$u_T(t):=u^+_T(t)-u^-_T(t)$,
where  $u^+_T, u^-_T\in \mathbb{M}\cap \mathbb{K}$ for all $T$.
Given $f\in\mathbb{M}$, set
$$\mathbb{B}_f:=\{g\in\mathbb{L}: g(t_2)-g(t_1)\geq f(t_2)-f(t_1), \mbox{ for all }
\ 0\leq t_1<t_2\leq 1\}.$$
The following lemma holds true.

\begin{lemma} \label{l4.4}
Suppose that for  $f\in \mathbb{V}$ we have
\beq \label{11.02.1}
\lim\limits_{T\rightarrow\infty}\int\limits_0^1|u_T(t)-f(t)|dt=0.\eeq
Then, for functions $f^\pm\in\mathbb{M}$ figuring in decomposition \eqref{fpm},
$$
\lim\limits_{T\rightarrow\infty}\inf_{g\in\mathbb{B}_{f^+}}\int\limits_0^1|u^+_T(t)-g(t)|dt=0, \ \ \
\lim\limits_{T\rightarrow\infty}\inf_{g\in\mathbb{B}_{f^-}}\int\limits_0^1|u^-_T(t)-g(t)|dt=0.$$
\end{lemma}

\noindent
{\sc Proof.} Let us prove that
$$\lim\limits_{T\rightarrow\infty}\inf_{g\in\mathbb{B}_{f^+}}\int\limits_0^1|u^+_T(t)-g(t)|dt=0.$$
Suppose the opposite, then there exists $\gamma>0$ such that for any $M>0$ there
exists $T>M$ such that
\beq \label{11.02.2}
\inf_{g\in\mathbb{B}_{f^+}}\int\limits_0^1|u^+_T(t)-g(t)|dt\geq \gamma.\eeq
Since functions $u^+_T$ lie in a compact, inequality \eqref{11.02.2} implies that there
exist a sub-sequence $T_M$ and a function $\tilde{g}$ such that
$$\lim\limits_{M\rightarrow\infty}\int\limits_0^1|u^+_{T_M}(t)-\tilde{g}(t)|dt=0, \ \ \
\inf_{g\in\mathbb{B}_{f_+}}\int\limits_0^1|\tilde{g}(t)-g(t)|dt \geq \gamma $$
and $\tilde{g}_\mathbb{D}\in \mathbb{M}$ because functions $u^+_{T_M}$,
$M=1,2,\ldots$, are monotone in $t$.

Therefore, from \eqref{11.02.1} it follows that
\beq\label{02.03.1}
\lim\limits_{M\rightarrow\infty}\int\limits_0^1|u^-_{T_M}(t)-(\tilde{g}(t)-f(t))|dt=0.\eeq
Then, since $u^-_T\in \mathbb{M}$, from \eqref{02.03.1} it follows that
$$\hat{g}_\mathbb{D}(t):=\tilde{g}_\mathbb{D}(t)-f_\mathbb{D}(t)$$
also belongs to $\mathbb{M}$. Hence, $f_\mathbb{D}(t)=\tilde{g}_\mathbb{D}(t)-\hat{g}_\mathbb{D}(t)$
where $\tilde{g}_\mathbb{D}\not\in B_{f^+}$
and $\hat{g}_\mathbb{D}\in\mathbb{M}$, which contradicts Lemma \ref{l4.3}.

In a similar fashion one can prove that
$$\lim\limits_{T\rightarrow\infty}\inf_{g\in\mathbb{B}_{f^-}}\int\limits_0^1|u^-_T(t)-g(t)|dt=0.
\quad \Box$$ \vspace{0.5cm}

A direct corollary of Lemma \ref{l4.4} is   \vskip 0.5cm

\begin{lemma} \label{c2}  Let $\mathbb{K}$ be a compact set in $(\mathbb{L},\rho )$. There exists
$\delta(\veps)>0$ such that $\lim\limits_{\veps\rightarrow 0}\delta(\veps)=0$
and for every $u\in\mathbb{K}\cap  \mathbb{U}_\veps(f)$ and $u^+$, $u^-$ from the decomposition
$u=u^+-u^-$ (cf.  \eqref{fpm}) the distances between $u^\pm$ and $\mathbb{B}_{f^\pm}$ satisfy
$$\rho(u^+,\mathbb{B}_{f^+})<\delta(\veps), \ \  \rho(u^-,\mathbb{B}_{f^-})<\delta(\veps).$$
\end{lemma}
\vspace{0.5cm}

\begin{lemma} \label{l4.6} Suppose function $u\in\mathbb{V}$ is increasing on $[0,1]$,
Let
$$\mathbb{B}_\veps:=\{g\in \mathbb{M}:\rho(g,u)<\veps\}.$$
Then there exists $\delta(\veps)>0$ such that
$$\inf\limits_{g\in\mathbb{B}_\veps}g(1)\geq u_\mathbb{D}(1)-\delta(\veps)$$
and \ $\lim\limits_{\veps\rightarrow 0}\delta(\veps)=0$.
\end{lemma}

\noindent
{\sc Proof.} Since $u_\mathbb{D}$ is increasing and left-continuous at $t=1$, there exists
a function $\gamma(\Delta)>0$ such that $\lim\limits_{\Delta\rightarrow 0}\gamma(\Delta)=0$
and
$$\sup\limits_{t\in[1-\gamma(\Delta),1]}(u_\mathbb{D}(1)-u_\mathbb{D}(t))<\Delta.$$
Let us choose $\Delta(\veps)$ so that $\Delta(\veps)\gamma(\Delta(\veps))\geq
\veps$ and $\lim\limits_{\veps\rightarrow 0}\Delta(\veps)=0$. Put
$$\delta(\veps):=3\Delta(\veps).$$
Suppose that
$$\inf\limits_{g\in\mathbb{B}_\veps}g(1)< u_\mathbb{D}(1)-\delta(\veps).$$
Then the condition $\inf\limits_{g\in\mathbb{B}_\veps}\rho(g,u)<\veps$ implies that there exists
a function $g\in\mathbb{B}_\veps$ such that
$$\begin{aligned}
\veps > & \diy\int\limits_{1-\gamma(\Delta(\veps))}^1|g(t)-u_\mathbb{D}(t)|dt\geq
\int\limits_{1-\gamma(\Delta(\veps))}^1
\bigl( |g(t)-u_\mathbb{D}(1)|-|u_\mathbb{D}(1)-u_\mathbb{D}(t)| \bigr)dt\\
> & \int\limits_{1-\gamma(\Delta(\veps))}^1|g(t)-u_\mathbb{D}(1)|dt-\Delta(\veps)\gamma(\Delta(\veps))\\
\geq & \int\limits_{1-\gamma(\veps)}^1|g(1)-u_\mathbb{D}(1)|dt-\Delta(\veps)\gamma(\Delta(\veps))
 >  2\gamma(\Delta(\veps))\Delta(\veps)-\Delta(\veps)\gamma(\Delta(\veps))>\veps.
\end{aligned}
$$
This contradiction completes the proof of the lemma. \quad $\Box$

\vspace{0.5cm}

\begin{lemma} \label{l4.5} {\rm{(\textrm{The ET property})}} Let condition
{\rm{(\ref{11.02.7})}} or {\rm{(\ref{17.12.6})}} be satisfied. Then for any $C>0$ there exists a
set $\mathbb{K}_C\subset\mathbb{L}$ compact in $(\mathbb{L},\rho )$ such that
$$\limsup_{T\rightarrow \infty}\dfrac{1}{\psi(T)}\ln\mathbf{P}\bigl( \xi_T(\,\cdot\,)\in
\mathbb{K}_C^\complement \bigr)\leq-C,$$
where $\mathbb{K}_C^\complement =\mathbb{L}\setminus\mathbb{K}_C$ and $\psi (T)=\varphi(T)\ln\varphi (T)$.
\end{lemma}

\noindent
{\sc Proof.} Take $\mathbb{K}_C:=\mathbb{V}_{a(C)}$ where $a(C) :=\dfrac{3C}{1-l}$. Then
\beq\label{11.02.4}\begin{aligned}
\mathbf{P}(\xi_T\in\mathbb{K}_C^\complement)&\leq
e^T\mathbf{E} \Bigl(e^{B_T+N_T\ln2}\mathbf{1}\Big[\zeta_T\in\mathbb{K}_C^\complement,
\inf\limits_{t\in[0,1]} \zeta_T(t)\geq 0\Big]\Bigr) \phantom{\sum_{r=\lfloor a(C)\varphi(T)\rfloor}}\\
&\leq e^T\mathbf{E}\Bigl(e^{B_T+N_T\ln2}\mathbf{1}\Big[N_T\geq a(C)\varphi(T), \inf\limits_{t\in[0,1]}
\zeta_T(t)\geq 0\Big]\Bigr)\\
& = e^T \sum\limits_{r=\lfloor a(C)\varphi(T)\rfloor}^\infty\mathbf{E}\Bigl(e^{B_T+N_T\ln2}
\mathbf{1}\Big[N_T= r, \inf\limits_{t\in[0,1]} \zeta_T(t)\geq 0\Big]\Bigr) \\
&\leq e^T \sum\limits_{r=\lfloor a(C)\varphi(T)\rfloor}^\infty \mathbf{E}
\Bigl(e^{B_T+N_T\ln2}\mathbf{1}\Big[N_T= r, \zeta^+(T)\geq\frac{r}{2}\Big]\Bigr),\end{aligned}\eeq
where the first inequality comes from \eqref{17.12.8} removing the $A_T$; the second inequality comes from
the observation that the process should have at least $a(C)\varphi(T)$ jumps during time interval $[0,T]$
to belong to the set $\mathbb{K}_C^\complement$. The last inequality means, that if the number of jumps
in the time interval $[0,1]$ is $r$, then to guarantee the inequality $\inf\limits_{t\in[0,1]} \zeta_T(t)\geq 0$ the
number of positive jumps should be at least $r/2$.

Let us upper-bound $B_T$ on the event
$\left\{\omega: N_T= r, \zeta^+(T)\geq\dfrac{r}{2}\right\}$ with $r\geq \lfloor a(C)\varphi(T)\rfloor$.
From condition \eqref{17.12.3} it follows that for any $\gamma>0$ and $T$ large enough,
\beq\label{11.02.3} \begin{aligned}
B_T=& \sum\limits_{i=1}^{r}\ln\bigl(\nu(\zeta(t_{i-1}),\zeta(t_i)\bigr) \\
\leq &\
 \zeta^-(T)\max\limits_{1 \leq i \leq r}\ln (\mu(i)) +\zeta^+(T) \max\limits_{1 \leq i \leq r}\ln (\lambda (i)) \phantom{\frac{r}{2}}\\
\leq & \ \zeta^-(T)\ln\bigl((1+\gamma)Q r\bigr) +\zeta^+(T) \ln\bigl((1+\gamma)P r^l\bigr) \phantom{\frac{r}{2}} \\
= & \ \zeta^-(T) \ln r + l\zeta^+(T) \ln r + r \ln M
\phantom{\frac{r}{2}}\\
= & \ \bigl( r - (1-l) \zeta^+(T) \bigr) \ln r + r \ln M
\leq \frac{r}{2}(1+l)\ln r+ r \ln M.\end{aligned}\eeq
Here $M:= (1+\gamma)^2(Q\vee 1)(P\vee 1)$.

By using \eqref{11.02.4}, \eqref{11.02.3} and the Stirling formula, we obtain that for
$r\geq \lfloor a(C)\varphi(T)\rfloor$ and $T$ large enough
\beq\label{11.02.5}\begin{aligned}
\mathbf{E}\Bigl( & e^{B_T+N_T\ln2}\mathbf{1}\Big[N_T= r, \zeta^+(T)\geq\dfrac{r}{2}\Big]\Bigr)\\
& \leq\exp\left\{\dfrac{r}{2}(1+l)\ln r+ r \ln (2M)\right\}\mathbf{P}(N_T= r)\\
& \leq e^{-T}\exp\left\{\dfrac{r}{2}(1+l)\ln r-r\ln r+ r \ln (2TMe)\right\}\\
& = e^{-T}\exp\left\{-\dfrac{r}{2}(1-l)\ln r+r \ln (2TMe)\right\}
\leq e^{-T} \exp\left\{-\dfrac{r}{3}(1-l)\ln r\right\},\end{aligned}
\eeq
where the last inequality is a consequence of the fact that under any of the conditions (\ref{11.02.7}) or (\ref{17.12.6}) the term $r\ln(2TMe)$ is $o(\frac r2(1-l)\ln(r))$ as $T$ tends to infinity.
The inequalities \eqref{11.02.4}, \eqref{11.02.5} imply that
$$
\begin{aligned} &\limsup\limits_{T\rightarrow \infty}\frac{1}{\psi(T)}\ln\mathbf{P}(\xi_T\in \mathbb{K}_C^\complement)\\
& \leq
\lim\limits_{T\rightarrow \infty}\frac{1}{\varphi(T)\ln\varphi(T)}
\ln\exp\Bigl\{-\dfrac{\lfloor a(C)\varphi(T)\rfloor }{3}(1-l)\ln \varphi(T)\Bigr\}=-C.\quad \Box\end{aligned}$$

\vspace{0.5cm}

Set
$$g_k:=e^{k\ln\varphi(T)}\dfrac{e^{-T/2}(T/2)^k}{k!}.$$

\begin{lemma} \label{l5.7} For any $C>0$ and $T>2C$
$$\max\limits_{0\leq k\leq C\varphi(T)}g_k= g_{\lfloor C\varphi(T)\rfloor }.$$
\end{lemma}

\noindent
{\sc Proof.} Given $1\leq k\leq \lfloor C\varphi(T) \rfloor$ where $T>2C$, we have
$$\dfrac{g_{k}}{g_{k-1}}=e^{\ln\varphi(T)}\dfrac{T/2}{k}=\dfrac{\varphi(T)T}{2k}\geq
\dfrac{\varphi(T)T}{2 \lfloor C\varphi(T)\rfloor }>1.$$
Thus, the sequence $g_k$ increases for $0\leq k\leq \lfloor C\varphi(T)\rfloor$. \quad $\Box$
\vskip .5cm

{\bf Acknowledgement.} The authors thank the referees for critical remarks and suggestions.
AL and AY thanks FAPESP for support under Grant 2022/01030-0 and 2017/10555-0. AL thanks IME, Universidade de Sao Paulo, for hospitality. YS thanks Math Department, Penn State University, for hospitality and support.
YS thanks St John's College, Cambridge, for support.

\end{document}